\documentclass[12pt]{article}
\usepackage{amstext}
\usepackage{amssymb}
\usepackage{amsfonts}
\usepackage{amsthm}
\usepackage{amsopn}
\usepackage{amsbsy}
\usepackage{layout}

\def\newtheorems{\newtheorem{theorem}{Theorem}[section]
                 \newtheorem{cor}[theorem]{Corollary}
                 \newtheorem{prop}[theorem]{Proposition}
                 \newtheorem{lemma}[theorem]{Lemma}
                 \newtheorem{defn}[theorem]{Definition}
                 
                 \newtheorem{notation}[theorem]{Notations}

                 \newtheorem{question}[theorem]{Question}
                 }

\newtheorems

\font\nmini=cmr10 scaled 550

\font\tinier=cmsy6 scaled 730

\font\srbf=cmbx8 scaled 920
\font\ssbf=cmbx8 scaled 700
\font\scmu=cmu10 scaled 600
\font\sscmu=cmu8 scaled 600

\font\blcmssifont=lcmssi8 scaled 1350

\def\fs#1{\mbox{\it #1\kern 1.3pt}}
\def\rfs#1{\mbox{\rm #1\kern 1.3pt}}
\def\bfs#1{\mbox{\bf #1\kern 1.3pt}}
\def\fss#1{\mbox{\scriptsize\it #1\kern 1.3pt}}
\def\fst#1{\mbox{\tiny\it #1\kern 1.1pt}}
\def\sifs#1{\mbox{\scriptsize\it #1\kern 1.3pt}}
\def\3norm#1{|\kern-1.3pt|\kern-1.3pt| #1 |\kern-1.3pt|\kern-1.3pt|}
\def\srfs#1{\mbox{\kern0.7pt\scriptsize\rm #1\kern 1.3pt}}
\def\sbfs#1{\mbox{\kern0.7pt\srbf #1\kern -0.6pt}}
\def\srbfs#1{\mbox{\kern0.7pt\srbf #1\kern -0.6pt}}
\def\spfs#1{\mbox{\kern0.7pt\scmu #1\kern 1.3pt}}
\def\sspfs#1{\mbox{\kern0.5pt\sscmu #1\kern 1.1pt}}
\def\ssbfs#1{\mbox{\kern0.7pt\ssbf #1\kern 1.3pt}}
\def\fsm#1{\mbox{\tiny\it #1\kern 1.0pt}}

\newcommand{\comp}{\hbox{$<\kern -3pt >$}}
\newcommand{\ncomp}
        {\;\hbox{\hbox{/}\kern -9.5pt \hbox{$<\kern -3pt >$}}}

\newcommand{\meet}
           {\hbox{$\wedge \kern -5.75pt \raise 1.5pt \hbox{$.$}\,$}}
\newcommand{\Meet}
         {\hbox{$\bigwedge \kern -8pt \raise 0.75pt \hbox{$.$}\:$}}
\newcommand{\ld}
           {\hbox{$< \kern -6pt \raise 2pt \hbox{$.$}\,$}}
\newcommand{\sss}{\: \hbox{$
\underline{\hbox{$\subset$}}\kern -4pt\raise -2pt \hbox{$\tiny |$}
$}\: }
\newcommand{\almostcontained}{\: \hbox{$
\raise 1.5pt \hbox{\scriptsize $\subset$}\kern -6.3pt\raise -3.5pt
\hbox{\scriptsize $\sim$}
$}\: }
\newcommand{\rraro}[2]{\hbox{$\kern 3pt\raise 2pt \hbox{$\raro$}
 \kern -14pt \raise
-3.5pt\hbox{\tiny{$#1\raro #2$}}$}}

\newcommand{\ct}{\centerline}

\newcommand{\frc}{\hbox{$\parallel \kern -5.7pt \hbox{$-$}$}}

\newcommand{\mrestriction}{\kern-1.1pt\upharpoonright\kern-1.1pt}
\newcommand{\nrestriction}{\kern-2.5pt\upharpoonright\kern-2.5pt}
\newcommand{\nnrestriction}{\kern-4.5pt\upharpoonright\kern-4.5pt}

\newcommand{\nfrc}{\not \kern -5pt \frc}
\newcommand{\rest}{\vbox{\hbox{$\:\kern -2pt\mathbin{\vert\kern-3.1pt\lower-1pt
   \hbox{$\mathsurround=0pt\mathchar"0012$}\kern-4pt}\:$}}}

\newcommand{\drest}{\vbox{\hbox{$\:\kern -2pt\mathbin{\rest\kern-4.6pt\lower1.7pt
   \hbox{$\mathsurround=0pt\mathchar"0012$}\kern-4pt}\:$}}}
\newcommand{\vdrest}{\vbox{\hbox{$\:\kern -2pt\mathbin{\rest\kern-3.95pt\lower1.7pt
   \hbox{$\mathsurround=0pt\mathchar"0012$}\kern-4pt}\:$}}}
   \newcommand{\smdrest}{\vbox{\hbox{$\:\kern -2pt\mathbin{\smrest\kern-2.7pt\lower1.7pt
      \hbox{$\mathsurround=0pt\mathchar"0012$}\kern-4pt}\:$\kern 1pt}}}

\newcommand{\rests}{\vbox{\hbox{\scriptsize$\:\kern
-1.4pt\mathbin{\vert\kern-2.4pt\lower-1pt
\hbox{\scriptsize$\mathsurround=0pt\mathchar"0012$}\kern-
3.0pt}\:$}}}

\newcommand{\pcntda}{\lower5pt\hbox{$\stackrel{\subset}{\neq}$}}
\newcommand{\pcntdb}{\lower5pt\hbox{$\stackrel{\supset}{\neq}$}}
\newcommand{\pcntdc}{\lower5pt\hbox{$\stackrel{\subseteq}{\neq}$}}

\newcommand{\m}[1]{\hbox{$ #1 $}}

\newcommand{\itGamma}{{\it \Gamma}}

\newcommand{\gammaprime}{\gamma\kern1pt'}

\newcommand{\calA}{{\cal A}}

\newcommand{\calB}{{\cal B}}

\newcommand{\calF}{{\cal F}}

\newcommand{\calN}{{\cal N}}

\newcommand{\calS}{{\cal S}}

\newcommand{\calU}{{\cal U}}

\newcommand{\calV}{{\cal V}}

\newcommand{\calW}{{\cal W}}

\newcommand{\hath}{\hat{h}}

\newcommand{\whatG}{\widehat{G}}

\newcommand{\whatP}{\widehat{P}}

\newcommand{\whatS}{\widehat{S}}

\newcommand{\swhatX}{\kern6pt\widehat{\rule{0pt}{8pt}}\kern-5.9pt X}

\newcommand{\whatGamma}{\widehat{\it \Gamma}}

\newcommand{\swtildeX}{\kern6pt\widetilde{\rule{0pt}{8pt}}\kern-5.9pt X}

\newcommand{\raro}{\rightarrow}

\newcommand{\itmb}[1]{\item[[\kern 1.5pt #1\kern -4pt]]}
\newcommand{\itms}[1]{\item[[#1\kern -5pt]]}

\newcommand{\II}{{\bf I\kern -1pt I}}

\newcommand{\third}{\frac{1}{3}}

\newcommand{\eqdf}{\hbox{\bf \,:=\,}}

\newcommand{\surj}{\vbox{\hbox{$\longrightarrow $
                  \kern -22pt \hbox{\lower 2.5pt  \hbox{\tiny onto}}
                  \kern -16pt \hbox{\raise 5pt  \hbox{\tiny 1-1}}
                  \kern 3pt}}}
\newcommand{\uarrow}[2]{\vbox{\hbox{$\longrightarrow $
                  \kern -16pt \hbox{\raise 5pt  \hbox{\tiny $#1$}}
                  \kern 10pt }}}

\newcommand{\spri}{\vbox{\hbox{\raise 2pt \hbox{\tiny $\|$}}}}
\newcommand{\tlr}{\vbox{\hbox{\raise 2pt \hbox{\tiny $\leftarrow$}}}}
\newcommand{\trr}{\vbox{\hbox{\raise 2pt \hbox{\tiny $\rightarrow$}}}}
\newcommand{\spr}{\mathrel{\vbox{\hbox{\tlr \kern -1.7pt \spri \kern -1.7pt \trr}}}}

\newcommand{\mspri}{\vbox{\hbox{\raise 2pt \hbox{$\scriptscriptstyle \|$}}}}

\newcommand{\mtlr}{\vbox{\hbox{\raise 2pt \hbox{$\scriptscriptstyle \leftarrow$}}}}

\newcommand{\mtrr}{\vbox{\hbox{\raise 2pt \hbox{$\scriptscriptstyle \rightarrow$}}}}

\newcommand{\mspr}{\mathrel{\vbox{\hbox{\mtlr \kern -1.7pt \mspri \kern -1.7pt \mtrr}}}}

\newcommand{\Nm}{\hbox{\kern -1.3pt \em I\kern-.2300em N\,}}
\newcommand{\Na}{\hbox{\kern -1.3pt \it I\kern-.2300em N\,}}
\newcommand{\Aa}{\hbox{\it A\kern -6.8pt\lower 1.0pt\hbox{-}\,}}
\newcommand{\Ba}{\hbox{\kern -1.3pt \it I\kern-.2300em B\,}}
\newcommand{\Da}{\hbox{\kern -1.3pt \it I\kern-.2300em D\,}}
\newcommand{\Ka}{\hbox{\kern -1.3pt \it I\kern-.2300em K\,}}
\newcommand{\La}{\hbox{\kern -1.3pt \it I\kern-.2300em L\,}}
\newcommand{\Ta}{\hbox{\kern 0.5pt \it T\kern-.5550em T\,}}
\newcommand{\jTa}{\hbox{\kern 0.5pt \it T\kern-.5550em T\,}}
\newcommand{\nTa}{\hbox{\kern 0.5pt \it T\kern-.6300em T\,}}
\newcommand{\Nmt}{\hbox{\kern -1.3pt I\kern-.2300em N\,}}
\newcommand{\N}{\hbox{I\kern-.1500em \hbox{\sf N}}}
\newcommand{\As}{\hbox{\scriptsize\it A\kern -5.3pt\lower 0.7pt\hbox{\it -}\,}}
\newcommand{\Ns}{\hbox{{\scriptsize\it I\kern-.1500em N}}}
\newcommand{\Nss}{\hbox{{\tiny\it I\kern-.1500em N}}}
\newcommand{\Ls}{\hbox{{\scriptsize\it I\kern-.1700em L}}}
\newcommand{\Rm}{\hbox{\kern -1.3pt \em I\kern-.1950em R\,}}
\newcommand{\Ra}{\hbox{\kern -1.3pt \it I\kern-.1950em R\,}}
\newcommand{\R}{\hbox{I\kern-.1500em \hbox{\sf R}}}
\newcommand{\sR}{\hbox{\tiny \hbox{I\kern-.1500em \hbox{\sf R}}}}
\newcommand{\msR}{\hbox{\tiny \hbox{\it I\kern-.2200em \hbox{\it R}}}}
\newcommand{\Rs}{\hbox{\scriptsize\it \hbox{I\kern-.2100em \hbox{R}}}}

\newcommand{\Q}
   {\hbox{${\rm Q} \kern -7.5pt \raise 2pt \hbox{\tiny$|$}\kern 7.5pt$}}
\newcommand{\C}
   {\hbox{${\rm C} \kern -6.5pt \raise 2pt \hbox{\tiny$|$}\kern 6.5pt$}}
\newcommand{\sC}
   {\hbox{\tiny \hbox{${\rm C} \kern -5.0pt \raise 1pt \hbox{\tiny$|$}\kern 7.5pt$}}}
\newcommand{\Z}{\hbox{\sf Z\kern-0.720em\hbox{ Z}}}
\newcommand{\sZ}{\hbox{\tiny\hbox{ \sf Z\kern-0.720em\hbox{ Z}}}}

\newcommand{\myqed}{\kern 5pt\vrule height7.5pt width6.9pt depth0.2pt
\rule{1.3pt}{0pt}}

\newcommand{\proofend}
{
\hbox{
\rule{0.8pt}{7pt}
\kern-3.9pt
\raise6.3pt
\hbox{\rule{5.0pt}{0.8pt}}
\kern-3.8pt
\rule{0.8pt}{7pt}
\kern-9.8pt
\rule{5.8pt}{0.8pt}
}
\kern-06pt}

\newcommand{\proofends}{\proofend\kern-07pt}

\newcommand{\proofendeol}
{
\hbox{
\rule{0.8pt}{7pt}
\kern-3.9pt
\raise6.3pt
\hbox{\rule{4.1pt}{0.8pt}}
\kern-3.8pt
\rule{0.8pt}{7pt}
\kern-9.0pt
\rule{4.9pt}{0.8pt}
\kern-5pt
}
\kern-8pt
}

\newcommand{\bd}{\begin{description}}
\newcommand{\ed}{\end{description}}

\newcommand{\ben}{\begin{enumerate}}
\newcommand{\een}{\end{enumerate}}

\newcommand{\sngltn}[1]{\{ #1 \}}

\newcommand{\setm}[2]{\{#1\:|\:#2\}}

\newcommand{\fsetn}[2]{\{\,#1,\ldots ,#2\}}

\newcommand{\norm}[1]{\| #1 \|}
\newcommand{\snorm}[1]{\hbox{\tiny$\|$} #1 \hbox{\tiny$\|$}}
\newcommand{\lnorm}[1]{\lower2pt\hbox{\scriptsize$\norm{#1}$}}
\newcommand{\boldnorm}[1]
{
\kern1pt\lower3pt\hbox{\rule{1.15pt}{12pt}\kern1.4pt}
#1
\lower3pt\hbox{\kern1.6pt\rule{1.15pt}{12pt}\kern2.0pt}
                                            }

\newcommand{\bboldnorm}[1]
{
\kern1pt\lower3pt\hbox{\rule{1.2pt}{12pt}\kern1.4pt}
#1
\lower3pt\hbox{\kern1.6pt\rule{1.2pt}{12pt}\kern2.0pt}
                                            }

\newcommand{\sboldnorm}[1]
{
\kern1pt\lower1.6pt\hbox{\rule{0.9pt}{07pt}\kern1.1pt}
#1
\lower1.6pt\hbox{\kern1.1pt\rule{0.9pt}{07pt}\kern1.8pt}
                                            }

\newcommand{\sbboldnorm}[1]
{
\kern1pt\lower0.5pt\hbox{\rule{1.0pt}{06pt}\kern1.1pt}
#1
\lower0.5pt\hbox{\kern1.1pt\rule{1.0pt}{06pt}\kern1.8pt}
                                            }
\newcommand{\abs}[1]{| #1 |}

\newcommand{\pair}[2]{\langle#1 ,#2\kern1.4pt\rangle}

\newcommand{\rpair}[2]{(#1 ,#2)}
\newcommand{\rfourtpl}[4]{(#1,#2,#3,#4)}

\newcommand{\fortpl}[4]{\langle#1 ,#2 ,#3 ,#4\kern1.4pt \rangle}

\newcommand{\change}[1]{\lower 0.7pt \hbox{\mbox{\large $#1$}}}

\newcommand{\fnn}[3]{#1:#2 \raro #3}
\newcommand{\tendin}[3]
{\lim_{#1\kern 1pt \ni\kern 1.5pt #2\kern 1pt \rightarrow\kern 1.5pt #3}}
\newcommand{\tend}[2]
{\lim_{#1\kern 1pt \rightarrow\kern 1.5pt #2}}

\newcommand{\iso}[3]{\m{ #1 : #2 \cong #3 }}

\newcommand{\authors}[3]{\bigskip\ct {#1} \bigskip\ct{#2}
\bigskip\ct{#3} \bigskip}

\newcommand{\inverse}{^{-1}}

\newcommand{\tsubseteq}{\lower2.1pt\hbox{\tiny$\subseteq$}}
\newcommand{\temptyset}{\lower1.6pt\hbox{\tiny$\emptyset$}}
\newcommand{\tSml}{\lower1.6pt\hbox{\tiny\it Sml}}
\newcommand{\tSpprtd}{\lower1.6pt\hbox{\tiny\it Spprtd}}
\newcommand{\tprec}{\lower1.6pt\hbox{\tiny$\prec$}}
\newcommand{\tcong}{\lower1.6pt\hbox{\tiny$\cong$}}
\newcommand{\tspr}{\lower3.5pt\hbox{\tiny$\spr$}}

\newcommand{\onetoone}{\hbox{1 \kern-3.2pt -- \kern-3.0pt 1}}

\newcommand{\onetoonen}
{\hbox
{1 \kern-3.2pt \rule[3.7pt]{4pt}{0.6pt} \kern-2.8pt 1}}

\newcommand{\num}[1]{$(#1)$}

\newcommand{\barg}{\bar{g}}
\newcommand{\barh}{\bar{h}}

\newcommand{\barE}{\bar{E}}

\newcommand{\overbbA}
               {\kern3.1pt
           \overline{\kern-1.1pt\bbA\kern-0.6pt}\kern0.6pt}
\newcommand{\oversbbA}
               {\kern2.6pt
           \overline{\kern-2.6ptsboldbbA\kern-0.3pt}\kern0.3pt}
\newcommand{\overB}
               {\kern3.1pt\overline{\kern-3.1ptB\kern-0.6pt}\kern0.6pt}
\newcommand{\dashB}
               {\kern3pt
           \overline{\phantom{\rule{1.0pt}{8.5pt}}}\kern1.2pt
               \overline{\phantom{\rule{1.0pt}{8.5pt}}}\kern1.2pt
               \overline{\phantom{\rule{1.0pt}{8.5pt}}}
           \kern-8.1pt B}
\newcommand{\overE}
               {\kern3.1pt\overline{\kern-3.1ptE\kern-0.6pt}\kern0.6pt}
\newcommand{\oversE}
               {\kern2.6pt\overline{\kern-2.6ptE\kern-0.3pt}\kern0.3pt}
\newcommand{\oversY}
               {\kern2.6pt\overline{\kern-1.7ptY\kern-0.3pt}\kern0.3pt}
\newcommand{\overssE}
               {\kern2.1pt
           \overline{\kern-1.2ptE\rule{0.3pt}{0pt}}
           \kern-0.3pt}
\newcommand{\overF}
               {\kern3.1pt\overline{\kern-3.1ptF\kern-0.6pt}\kern0.6pt}
\newcommand{\oversF}
               {\kern2.6pt\overline{\kern-2.6ptF\kern-0.3pt}\kern0.3pt}
\newcommand{\overG}
               {\kern3.1pt\overline{\kern-3.1ptG\kern-0.6pt}\kern0.6pt}
\newcommand{\overH}
               {\kern3.1pt\overline{\kern-3.1ptH\kern-0.6pt}\kern0.6pt}
\newcommand{\overL}
               {\kern3.1pt\overline{\rule{2pt}{0pt}\kern-3.1ptL\kern-0.6pt}\kern0.6pt}
\newcommand{\overP}
               {\kern3.1pt\overline{\kern-0.9ptP\kern-0.3pt}\kern0.6pt}
\newcommand{\overV}
               {\kern3.1pt\overline{\kern-0.9ptV\kern-0.3pt}\kern0.6pt}
\newcommand{\overW}
               {\kern3.1pt\overline{\kern-0.9ptW\kern-0.3pt}\kern0.6pt}
\newcommand{\overX}
               {\kern3.1pt\overline{\kern-3.1ptX\kern-0.6pt}\kern0.6pt}
\newcommand{\overY}
               {\kern3.1pt\overline{\rule{2pt}{0pt}\kern-3.1ptY\kern-0.6pt}\kern0.6pt}
\newcommand{\overZ}
               {\kern3.1pt\overline{\kern-3.1ptZ\kern-0.6pt}\kern0.6pt}
\newcommand{\oversG}
               {\kern2.6pt\overline{\kern-2.6ptG\kern-0.3pt}\kern0.3pt}
\newcommand{\oversX}
               {\kern2.6pt\overline{\kern-2.6ptX\kern-0.3pt}\kern0.3pt}
\newcommand{\oversZ}
               {\kern2.6pt\overline{\kern-2.6ptZ\kern-0.3pt}\kern0.3pt}

\newcommand{\overcalO}
               {\kern3.1pt\overline{\kern-3.1pt{\cal O}\kern-0.6pt}\kern0.6pt}

\newcommand{\sstar}{\hbox{\scriptsize \kern1pt$\star$\kern0pt}}
\newcommand{\rsstar}{\raise1pt\hbox{\scriptsize \kern1pt$\star$\kern1pt}}
\newcommand{\raisedsstar}{\kern0.5pt\raise1.0pt\hbox{\scriptsize$*$}}
\newcommand{\sperp}{\raise2pt\hbox{\kern0pt\tiny $\perp$}}
\newcommand{\ssperp}{\raise1.5pt\hbox{\kern0pt\tinier \symbol{63}}}

\newcommand{\nperp}
{
\hbox{\kern1.0pt
\hbox{\raise1.0pt \hbox{\rule{4.6pt}{0.3pt}}}\kern-5.1pt$\sperp$}
}
\newcommand{\snperp}
{
\hbox{\kern1.0pt
\hbox{\raise0.4pt \hbox{\rule{3.3pt}{0.2pt}}}\kern-3.65pt$\ssperp$}
}

\newcommand{\sor}{\raise2pt\hbox{\tiny\rm Or}}

\newcommand{\srs}[1]{\raise1.3pt\hbox{\tiny\rm #1}}
\newcommand{\ssrfs}[1]{\raise1.3pt\hbox{\tiny\rm #1}}
\newcommand{\ssrs}[1]{\raise1.3pt\hbox{\nmini #1}}
\newcommand{\scirc}
{\raise1pt\hbox{\scriptsize\kern1.5pt$\circ$\kern1.5pt}}
\newcommand{\sscirc}
{\hbox{\tiny\kern1.5pt$\circ$\kern0.3pt}}
\newcommand{\tcirc}
{\kern0.5pt\raise2pt\hbox{\tiny$\circ$}}

\newcommand{\bcirc}
{\mathop{\lower1pt\hbox{\large\kern1.5pt$\circ$}}}
\newcommand{\bbcirc}
{\mathop{\lower1.8pt\hbox{\Large\kern1.5pt$\circ$}}}

\newcommand{\raisedstar}
{\raise-3.1pt\hbox{$\rule{0.5pt}{0pt}^*$}}
\newcommand{\rstar}
{\raise-3.1pt\hbox{$\rule{0.5pt}{0pt}^*$}}
\newcommand{\raisedtinystar}
{\kern0.7pt\raise1.5pt\hbox{\tiny $*$}}
\newcommand{\rtstar}
{\kern0.7pt\raise1.0pt\hbox{\tiny $*$}\kern-2.0pt}
\newcommand{\bfie}{\hbox{\blcmssifont \symbol{101}\kern1pt}}
\newcommand{\bfif}{\hbox{\blcmssifont \symbol{102}\kern1pt}}

\newcommand{\bfih}{\hbox{\blcmssifont \symbol{104}\kern1pt}}

\newcommand{\mcdot}{\kern-1.4pt\cdot\kern-1.4pt}
\newcommand{\lcdot}{\hbox{$\kern0.1pt\cdot\kern1.0pt$}}
\newcommand{\ncdot}{\hbox{$\kern0.7pt\cdot\kern0.6pt$}}
\newcommand{\fr}{\kern1.0ptr}
\newcommand{\fzero}{\kern1.3pt0}
\newcommand{\fone}{\kern1.0pt1}
\newcommand{\ftwo}{\kern1.3pt2}
\newcommand{\fprime}{\kern1.3pt'}
\newcommand{\farsu}[1]{^{\kern1.3pt#1}}
\newcommand{\sprime}{\raise1.3pt\hbox{\scriptsize\kern1.4pt$'$}}
\newcommand{\fprimesub}{\kern1.5pt'\kern-3.8pt}
\newcommand{\fprimew}{\kern0.7pt'}
\newcommand{\fprimei}{\kern1.5pt'\kern-3.0pt}
\newcommand{\fnprime}{\kern1.5pt'\kern-3.0pt}

\newcommand{\clubsign}{{\raise1.3pt\hbox{\tiny $\clubsuit$}}}
\newcommand{\sclub}{{\raise1.3pt\hbox{\tiny $\clubsuit$}}}
\newcommand{\srsp}[1]{\hbox{\kern1pt\tiny $(#1)$}}
\newcommand{\sminus}{\rule{5.6pt}{0.3pt}}
\newcommand{\ssim}{\hbox{\scriptsize $\sim$}}
\newcommand{\neweq}{
\mathrel{
\hbox{
\lower-0.0pt\hbox{\sminus}
\kern-9.95pt\ssim\kern-6.2pt\raise3.8pt\hbox{\sminus}
}}}

\newcommand{\dcup}{\mathbin{\hbox{$\cup$
\kern-10.43pt\raise1.4pt \hbox{\tiny $\cup$}}}}

\newcommand{\sprt}[1]
{\lower1.3pt\hbox{\kern1.5pt\rule{0.5pt}{10pt}}
\kern-0.5pt\underline{\kern1pt#1\kern1pt}
\lower1.3pt\hbox{\kern0.2pt\rule{0.5pt}{10pt}}\kern1pt}

\newcommand{\sprtd}[2]
{#1\lower1.3pt\hbox{\kern1.5pt\rule{0.5pt}{10pt}}
\kern-0.5pt\underline{\kern1pt#2\kern1pt}
\lower1.3pt\hbox{\kern0.2pt\rule{0.5pt}{10pt}}\kern1pt}

\newcommand{\sprtl}[1]
{\lower4.4pt\hbox{\kern1.5pt\rule{0.5pt}{14pt}}
\kern-0.5pt\underline{\kern1pt#1\kern1pt}
\lower4.4pt\hbox{\kern-0.5pt\rule{0.5pt}{14pt}}\kern1pt}

\newcommand{\sprtll}[1]
{\lower5.8pt\hbox{\kern1.5pt\rule{0.5pt}{15.4pt}}
\kern-0.5pt\underline{\kern1pt#1\kern1pt}
\lower5.8pt\hbox{\kern-0.5pt\rule{0.5pt}{15.4pt}}\kern1pt}

\newcommand{\sprtdl}[2]
{#1\lower4.0pt\hbox{\kern1.5pt\rule{0.5pt}{14pt}}
\kern-0.5pt\underline{\kern1pt#2\kern1pt}
\lower4.1pt\hbox{\kern0.0pt\rule{0.5pt}{14pt}}\kern1pt}

\newcommand{\sprtm}[1]
{\lower3.3pt\hbox{\kern1.5pt\rule{0.5pt}{14pt}}
\kern-0.5pt\underline{\kern1pt#1\kern1pt}
\lower3.4pt\hbox{\kern0.0pt\rule{0.5pt}{14pt}}\kern1pt}

\newcommand{\sprtdm}[2]
{#1\lower3.3pt\hbox{\kern1.5pt\rule{0.5pt}{14pt}}
\kern-0.5pt\underline{\kern1pt#2\kern1pt}
\lower3.4pt\hbox{\kern0.0pt\rule{0.5pt}{14pt}}\kern1pt}

\DeclareSymbolFont{AMSb}{U}{msb}{m}{n}
\DeclareMathSymbol{\bbA}{\mathbin}{AMSb}{"41}
\DeclareMathSymbol{\bbB}{\mathbin}{AMSb}{"42}
\DeclareMathSymbol{\bbC}{\mathbin}{AMSb}{"43}
\DeclareMathSymbol{\bbD}{\mathbin}{AMSb}{"44}
\DeclareMathSymbol{\bbE}{\mathbin}{AMSb}{"45}
\DeclareMathSymbol{\bbF}{\mathbin}{AMSb}{"46}
\DeclareMathSymbol{\bbG}{\mathbin}{AMSb}{"47}
\DeclareMathSymbol{\bbH}{\mathbin}{AMSb}{"48}
\DeclareMathSymbol{\bbI}{\mathbin}{AMSb}{"49}
\DeclareMathSymbol{\bbJ}{\mathbin}{AMSb}{"4A}
\DeclareMathSymbol{\bbK}{\mathbin}{AMSb}{"4B}
\DeclareMathSymbol{\bbL}{\mathbin}{AMSb}{"4C}
\DeclareMathSymbol{\bbM}{\mathbin}{AMSb}{"4D}
\DeclareMathSymbol{\bbN}{\mathbin}{AMSb}{"4E}
\DeclareMathSymbol{\bbO}{\mathbin}{AMSb}{"4F}
\DeclareMathSymbol{\bbP}{\mathbin}{AMSb}{"50}
\DeclareMathSymbol{\bbQ}{\mathbin}{AMSb}{"51}
\DeclareMathSymbol{\bbR}{\mathbin}{AMSb}{"52}
\DeclareMathSymbol{\bbS}{\mathbin}{AMSb}{"53}
\DeclareMathSymbol{\bbT}{\mathbin}{AMSb}{"54}
\DeclareMathSymbol{\bbU}{\mathbin}{AMSb}{"55}
\DeclareMathSymbol{\bbV}{\mathbin}{AMSb}{"56}
\DeclareMathSymbol{\bbW}{\mathbin}{AMSb}{"57}
\DeclareMathSymbol{\bbX}{\mathbin}{AMSb}{"58}
\DeclareMathSymbol{\bbY}{\mathbin}{AMSb}{"59}
\DeclareMathSymbol{\bbZ}{\mathbin}{AMSb}{"5A}
\newcommand{\boldbbA}{\mathbb{A}\kern-9pt\mathbb{A}}
\newcommand{\boldbbL}{\mathbb{L}\kern-8.3pt\mathbb{L}}
\newcommand{\boldbbR}{\mathbb{R}\kern-8.3pt\mathbb{R}}

\newcommand{\sboldbbA}
                {\hbox{\kern0.7pt\scriptsize
        $\mathbb{A}\kern-6pt\mathbb{A}$}}
\newcommand{\sboldbbL}
                {\hbox{\kern0.7pt\scriptsize
        $\mathbb{L}\kern-5.68pt\mathbb{L}$}}
\newcommand{\sboldbbN}
                {\hbox{\kern0.7pt\scriptsize
        $\mathbb{N}\kern-5.60pt\mathbb{N}$}}
\newcommand{\sboldbbR}
                {\hbox{\kern0.7pt\scriptsize
        $\mathbb{R}\kern-5.68pt\mathbb{R}$}}
\newcommand{\sboldbbT}
                {\hbox{\kern0.7pt\scriptsize
        $\mathbb{T}\kern-5.60pt\mathbb{T}$}}

\newcommand{\wpm}{\kern0.7pt\pm}
\newcommand{\doubledagger}{\hbox{$\dagger\kern-4pt\dagger$}}

\newcommand{\Wr}{{\kern3pt\bf Wr\kern3pt}}

\begin{document}
\baselineskip 18pt
\setcounter{page}{1}
\centerline
{\Large\bf A\kern-2pt\ Reconstruction\kern-2pt\ theorem\kern-2pt\
for\kern-2pt\ homeomorphism}
\vspace{3mm}
\centerline
{\Large\bf groups without small sets and}
\vspace{3mm}
\centerline
{\Large\bf non-shrinking\kern-2pt\ functions of a normed space
}
\vspace{3mm}

    \bigskip

    \authors
{\bf By Vladimir P. Fonf and Matatyahu Rubin\vspace{-3mm}}
{Department of Mathematics Ben Gurion University\vspace{-3mm}}
{Beer Sheva, Israel}

\section{Introduction}

Let $X$ be a topological space and $G$ be a subgroup of
the group $H(X)$ of all auto-homeomorphisms of $X$.
The pair $\rpair{X}{G}$ is then called a {\it space-group pair}.
Let $K$ be a class of space-group pairs.
$K$ is called a {\it faithfull class} if for every
$\rpair{X}{G}, \rpair{Y}{H} \in K$
and an isomorphism $\varphi$ between the groups
$G$ and $H$ there is a homeomorphism $\tau$ between $X$ and $Y$
such that $\varphi(g) = \tau \scirc g \scirc \tau\inverse$
for every $g \in G$.

The first important theorem on faithfulness is due to J. Whittaker
\cite{W} (1963). He proved that the class of homeomorphism groups
of Euclidean manifolds is faithful.
That is, $\setm{\rpair{X}{H(X)}}{X \mbox{ is a Euclidean manifold}}$
is faithful.

Other faithfulness theorems were proved in
R. McCoy \cite{McC} 1972, W. Ling \cite{Lg1} 1980,
M. Rubin \cite{Ru1} 1989,
M. Rubin \cite{Ru2} 1989,
K. Kawamura \cite{Ka} 1995,
M. Brin \cite{Br1} 1996,
M. Rubin \cite{Ru3} 1996,
A. Banyaga \cite{Ba1} 1997,
A. Leiderman and R. Rubin \cite{LR} 1999
and
J. Borzellino and V. Brunsden 2000.

\newpage
Among the classes shown to be faithful in \cite{Ru2}
is the class of manifolds over normed spaces.
In \cite{RY} 2000, M. Rubin and Y. Yomdin
obtained various strengthenings and continuations of the
this result from \cite{Ru2}.
One central result from \cite{RY} is the following theorem.

{\bf Theorem A} Let $K_{\fss{LLIP}}$ be the class of all
space-group pairs $\rpair{X}{G}$ such that $X$ is an open subset of a
normed space and $G$ is a subgroup of $H(X)$ which contains
all locally bilipschitz homeomorphisms of $X$.
\newline
Then $K_{\fss{LLIP}}$ is faithful.

\kern2mm

Let $\rpair{X}{G}$ be a space-group pair
and $\emptyset \neq U \subseteq X$ be open.
We say that $U$ is a {\it small set}
with respect to $\rpair{X}{G}$,
if for every open nonempty $V \subseteq U$ there is $g \in G$ such that
$g(U) \subseteq V$.

It is easy to show that if $X$ is an open subset of a normed
space,
and $\fs{LIP}(X)$ is the group of all bilipschitz homeomorphisms
of $X$, then
the family of subsets of $X$ which are small with respect to
$\rpair{X}{\fs{LIP}(X)}$ is an open cover of $X$.
The same is obviously true for any group of homeomorphisms containing
$\fs{LIP}(X)$.
The existence of a cover consisting of small sets is indeed used in
the proof of Theorem~A. This fact is also used
in all previous faithfulness results applicable to
infinite dimensional normed spaces.

This leads to the question of discovering subgroups of $H(X)$
which are rich enough to allow the recovery of $X$, but which
are not sufficiently big to admit small sets.

This work addresses this question.
We prove a new faithfulness result (Theorem \ref{t1.2})
which does not assume the existence of small sets.
And we also construct a large class of groups which
do not have small sets, and which are covered by this new
faithfulness result.

Theorem \ref{t1.2} deals with a general class of first countable spaces.
We apply it to the class of open subsets of normed spaces.
It is also applicable outside the class of normed spaces.
One such application - to the class of metrizable locally
convex spaces, is proved in Theorem \ref{t4.9}.

The following statement is a special case of Theorem \ref{t1.2}.
Let $\rpair{X}{G}$ be a space-group pair and
$S \subseteq X$ be open. $S$ is {\it strongly flexible},
if for every infinite $A \subseteq S$ without
accumulation points in $X$,
there is a nonempty open set $V \subseteq X$ such that
for every nonempty open set $W \subseteq V$ there is $g \in G$
such that the sets
$\setm{a \in A}{g(a) \in W}$ and
$\setm{a \in A}
{\mbox{for some neighborhood } U \mbox{ of } a,\break
g \nrestriction U  = \fs{Id}} $
are infinite.

\kern0.7mm

{\bf Theorem B }
Let $K_F$ be the class of all space-group pairs $\rpair{X}{G}$
such that
\newline
(1) $X$ is regular, first countable and has no isolated points.
\newline
(2) For every $x \in X$ and an open neighborhood $U$ of $x$
the set
\newline\indent
$\setm{g(x)}{g \in G \mbox{ and }
g \nrestriction (X - U) = \fs{Id}}$
is somewhere dense.
\newline
(3) The family of strongly flexible sets is a cover of $X$.
\\
\underline{Then} $K_F$ is faithful.

\kern0.7mm

\kern4mm

\noindent
{\bf Normed spaces}

We next describe our results for normed spaces.
Let $\fs{LLIP}(X)$ denote the group of locally bilipschitz
homeomorphisms of a metric space $X$.
For every nonempty open subset $X$ of an infinite dimensional
normed space we shall define a certain subgroup
$G_X \subseteq \fs{LLIP}(X)$.
We shall prove the following theorem concerning the $G_X$'s.

\kern0.7mm

{\bf Theorem C }
(a) If $X$ is a nonempty open subset an infinite dimensional
normed space, then $X$ does not have small sets with respect to
$\rpair{X}{G_X}$.

(b) Let $K_{\fss{NSML}}$ be the class of all space-group pairs
$\rpair{X}{H}$
such that
$G_X \subseteq H \subseteq H(X)$.
Then $K_{\fss{NSML}}$ is faithful.

\kern0.7mm

The group $G_X$ is defined in Definition \ref{d4.2}.
Part (a) of Theorem~B is a corollary of Theorem \ref{t3.1}(a).
Part (b) of Theorem~B follows from Theorems \ref{t1.2} and \ref{t4.3}.

Note that $ K_{\fss{LLIP}} \subseteq K_{\fss{NSML}}$.
Recall that for every $\rpair{X}{H} \in K_{\fss{LLIP}}$,
the family of small sets is a cover of $X$.
On the other hand, $\rpair{X}{G_X}$ has no small sets at all,
and it belongs to $K_{\fss{NSML}}$.
So $K_{\fss{LLIP}} \subsetneqq K_{\fss{NSML}}$.
Hence Theorem~C(b) strengthens Theorem A.

It needs to be mentioned that if $E$ is an infinite dimensional
normed space and $X \subseteq E$ is open and nonempty,
then $G_X$ is obtained from $G_E$ in the following way.
\newline\centerline{
$G_X =
\setm{g \nrestriction X}{g \in G_E \mbox{ and }
g \nrestriction (E \setminus X) = \fs{Id}}$.
}
So it suffices to show that $\rpair{E}{G_E}$ does not have small sets.

\kern2mm

The nonexistence of small sets follows from a stronger property
which is called here the non-shrinking property.
For a metric space $X$, $x \in X$ and $r > 0$ let
$B(x,r)$ denote the closed ball with center at $x$ and radius $r$.
Let $X$ be a metric space and $\itGamma$ be a semigroup of functions
from $X$ to $X$. The operation in $\itGamma$ is composition of
functions.
$\itGamma$ is {\it non-shrinking} if there are no $g \in \itGamma$,
$x \in X$ and $r > 0$ such that $g(B(x,r))$ is contained in a
finite union of closed balls with radius $< r$.

The final result of Section~\ref{s3} is Theorem \ref{t3.1}.
For a normed space $E$ we shall define two semigroups of functions from
$E$ to $E$:
$\itGamma_1(E)$ and $\itGamma_2(E)$.
$\itGamma_1(E)$ is a subsemigroup of $\itGamma_2(E)$.
Theorem \ref{t3.1}(a) states that for every infinite dimensional
normed space $E$, $\itGamma_1(E)$ is non-shrinking.

Theorem \ref{t3.1}(b) says that if $E$
has an infinite dimensional subspace $F$ such that $c_0$
(the space of real sequences converging to $0$),
is not isomorphically embeddable in the completion of $F$,
then $\itGamma_2(E)$ is non-shrinking.

Also observed in Section \ref{s3} is that $\itGamma_2(c_0)$
is {\it not} non-shrinking.

The group $G_E$ from Theorem C is contained in $\itGamma_1(E)$.
So the fact that $G_E$ and hence $G_X$ has no small sets follows
from Theorem \ref{t3.1}(a).


\kern2mm


\kern1.0mm

\noindent
{\bf Metrizable locally convex spaces}

Theorem D is another application of Theorem \ref{t1.2}.
It is an analogue of Theorem A to the class of
metrizable locally convex spaces.

{\bf Theorem D } Let $K_M$ be the class of all space-group
pairs $\rpair{X}{G}$
in which $X$ is an open subset of a locally convex metrizable
topological vector space, and $G$ is a group containing all locally
bi-uniformly continuous homeomorphisms of $X$.
Then $K_M$ is faithful.

\kern0.7mm

Theorem~C is restated in \ref{t4.9}(a).

\kern1.0mm

\noindent
{\bf Acknowledgements }
Lemma \ref{l3.5} which is used in the proof of Theo-\break
rem~\ref{t3.1}(a)
was found by Michael Levin.
We thank him for his kind permission to include it.
We had another proof that did not rely on Lemma \ref{l3.5}.
But the proof which uses Levin's Lemma is simpler.

We also thank Arkady Leiderman for his help in the proof of
Proposition~\ref{p4.12}.

\section{The reconstruction theorem}\label{s2}

\begin{notation}
\begin{rm}
Let $X$ be a topological space.
If $A \subseteq X$. Then the closure and the interior of $A$ in $X$
are denoted respectively by $\fs{cl}_X(A)$ and $\fs{int}_X(A)$.
Also, $\fs{acc}_X(A)$ denotes the set of accumulation points of
$A$ in $X$.
If $x \in X$ then $\fs{Nbr}_X(x)$ denotes the set of open
neighborhoods of $x$ in $X$.
The subscript $X$ is sometimes omitted.

\end{rm}
\end{notation}

\begin{defn}\label{d1.1}
\begin{rm}
Let $X$ be a topological space.

(a)
A subset $A \subseteq X$ is {\it discrete} if $\fs{acc}(A) = \emptyset$.

(b) Let $\rpair{X}{G}$ be a space-group pair.
The group $G$ is called a {\it locally moving} group of $X$,
if for every open nonempty set $U \subseteq X$
there is\break
$g \in G \setminus \sngltn{\fs{Id}}$ such that
$g \nrestriction (X \setminus U) = \fs{Id}$.

(c) Let $\rpair{X}{G}$ be a space-group pair,
$A \subseteq X$ be infinite
and $V \subseteq X$ be open and nonempty.
We say that $V$ {\it dissects} $A$, if for for every open nonempty
$W \subseteq V$ there is $g \in G$ such that
$\setm{a \in A}{g(a) \in W}$ is infinite,
and $\setm{a \in A}{\mbox{there is } S \in \fs{Nbr}(a)
\mbox{ such that } g \nrestriction S = \fs{Id}}$
is infinite.
$A$ is {\it dissectable} if there is $V$ such that $V$ dissects $A$.

(d) Let $\rpair{X}{G}$ be a space-group pair and $U \subseteq X$ be
open. $U$ is {\it flexible} if there is a dense subset $D \subseteq U$
such that every infinite discrete subset of $D$ is dissectable.

(e) Let $\rpair{X}{G}$ be a space-group pair.
\newline
(1)
The set $D(X,G)$ is defined as follows. A point $x$ is in $D(X,G)$ iff
$\setm{g(x)}{g \in G}$ is somewhere dense.
That is, if $\fs{int}(\fs{cl}(\setm{g(x)}{g \in G})) \neq \emptyset$.
\newline
(2)
Define the relation $\fs{DF}_{X,G}(x,y)$ as follows.
$\fs{DF}_{X,G}(x,y)$ holds if
the set\break
\centerline{
$\setm{g(x)}{g \in G \mbox{ and there is } U \in \fs{Nbr}(y)
\mbox{ such that } g \nrestriction U = \fs{Id}}$
}
is somewhere dense.

(f) Define the class $K_1$ of space-group pairs as follows.
Let $\rpair{X}{G}$ be a space-group pair.
$\rpair{X}{G} \in K_1$ if the following hold:
\newline
(P1) $X$ is regular and first countable.
\newline
(P2) $G$ is a locally moving group of $X$.
\newline
(P3) For every distinct $x,y \in X$,
$\fs{DF}_{X,G}(x,y)$ holds.
\newline
(P4) The set of flexible subsets of $X$ is a cover of $X$.

Define the class $K$ of space-group pairs as follows.
Let $\rpair{X}{G}$ be a space-group pair.
$\rpair{X}{G} \in K$ if the following hold:
\newline
(P1) $X$ is regular and first countable.
\newline
(P2) $G$ is a locally moving group of $X$.
\newline
(Q1) $D(X,G)$ is a dense subset of $X$.
\newline
{\thickmuskip=2mu \medmuskip=1mu \thinmuskip=1mu
(Q2) For every distinct $x,y \in D(X,G)$,
$\fs{DF}_{X,G}(x,y)$ holds or $\fs{DF}_{X,G}(y,x)$ holds.
}
(P4) The set of flexible subsets of $X$ is a cover of $X$.
\end{rm}
\end{defn}

Note that $K_1 \subseteq K$.
Note also that $D(X,G)$ is invariant under $G$.

\begin{theorem}\label{t1.2}
\num{a} $K_1$ is faithful.

\num{b} For every $\rpair{X_1}{G_1}, \rpair{X_2}{G_2} \in K$ and
$\iso{\varphi}{G_1}{G_2}$
there is
\newline
$\iso{\tau}{D(X_1,G_1)}{D(X_2,G_2)}$
such that $\tau$ induces $\varphi$.
That is, for every $g \in G$,
$\varphi(g) \nrestriction D(Y,H) =
\tau \scirc (g \nrestriction D(X,G)) \scirc \tau\inverse$.
\end{theorem}

\noindent
{\bf Remark } Part (b) of Theorem \ref{t1.2} implies Part (a).
In this work only Part~(a) is used. But there are concrete classes
of space-group pairs, which require the use of Part (b).
An example of such a class, is the class of all space-group pairs
$\rpair{X}{G}$, in which $X$ is the closure of an open subset of a
normed space, and $G$ is the group of all homeomorphisms of
$X$ which take the boundary of $X$ to itself.

We shall use a theorem from \cite{Ru3}.
It is quoted here as Theorem \ref{t2.3}.
The set of regular open subsets of $X$
is denoted by $\fs{Ro}(X)$.
Recall that $U$ is a regular open set if $\fs{int}(\fs{cl}(U)) = U$.
For $U,V \in \fs{Ro}(X)$ define\break
$U + V = \fs{int}(\fs{cl}(U \cup V))$ and $\sim U = \fs{int}(X - U)$.
Then $\rfourtpl{\fs{Ro}(X)}{+}{\cap}{\sim}$ is a complete
Boolean algebra, which we denote by $\bfs{Ro}(X)$.
The partial ordering of the Boolean algebra $\bfs{Ro}(X)$ is
$\subseteq$.
Note that every $g \in H(X)$ induces an automorphism of
$\bfs{Ro}(X)$ which we also denote by $g$.

\begin{theorem}\label{t2.3}
Let $\pair{X}{G}$ and $\pair{Y}{H}$ be space-group pairs.
Assume that $G$ and $H$ are locally moving groups of $X$ and $Y$
respectively.
Let $\iso{\varphi}{G}{H}$.
Then there is a unique
$\iso{\eta}{\bfs{Ro}(X)}{\bfs{Ro}(Y)}$
such that for every $g \in G$,
$\varphi(g) = \eta \scirc g \scirc \eta\inverse$.
\end{theorem}

\noindent
{\bf Proof }
See \cite{Ru3} Definition 1.2, Corollary 1.4 or Corollary 2.10
and Proposition 1.8.
\hfill\myqed

\begin{defn}\label{d3.3}
\begin{rm}
(a) Let $\calB$ be a set of sets. $\calB$ is called a
{\it pairwise disjoint family} if for every distinct $A,B \in \calB$,
$A \cap B = \emptyset$.
$\calB$ is called a {\it regular open family},
if $\calB \subseteq \fs{Ro}(X)$.
Let $\calA$ be an infinite set of subsets of $X$ and $x \in X$.
$\lim \calA = x$, if for every $S \in \fs{Nbr}(x)$,
$\setm{A \in \calA}{A \not\subseteq S}$ is finite.
$\calA$ is {\it convergent} if for some $x \in X$, $\lim \calA = x$.

(b) Let $\calA$ be a set of subsets of $X$.
$\fs{acc}(\calA)$ is defined as follows:
$x \in \fs{acc}(\calA)$ iff there is $B \subseteq \bigcup \calA$
such that for every $A \in \calA$, $\abs{B \cap A} \leq 1$
and $x \in \fs{acc}(B)$.

(c) Let $\calU$ be an infinite pairwise disjoint regular open family
and $V \in \fs{Ro}(X) \setminus \sngltn{\emptyset}$.
$V$ {\it absorbs} $\calU$,
if for every $W \in \fs{Ro}(X) \setminus \sngltn{\emptyset}$:
if $W \subseteq V$, then there is $g \in G$ such that
$\setm{U \in \calU}{g(U) \not\subseteq W}$ is finite.

$\calU$ is {\it absorbable} if there is
$V \in \fs{Ro}(X) \setminus \sngltn{\emptyset}$
such that $V$ absorbs $\calU$.

(d) Let $\calU$ be an infinite pairwise disjoint regular open family
and $V \in \fs{Ro}(X) \setminus \sngltn{\emptyset}$.
$V$ {\it splits} $\calU$
if for every $W \in \fs{Ro}(X) \setminus \sngltn{\emptyset}$:
if $W \subseteq V$, then there is $g \in G$
such that
$\setm{U \in \calU}{g(U) \cap W \neq \emptyset}$ is infinite and
\newline
$\setm{U \in \calU}
{\mbox{there is } U' \in \fs{Ro}(X) \setminus \sngltn{\emptyset}
\mbox{ such that } U' \subseteq U
\mbox{ and } g \restriction U' = \fs{Id}}$\break
is infinite.

$\calU$ is {\it splittable} if there is
$V \in \fs{Ro}(X) \setminus \sngltn{\emptyset}$
such that $V$ splits $\calU$.
\end{rm}
\end{defn}

\begin{prop}\label{p1.3}
Let $\calU$ be an infinite pairwise disjoint regular open family
and $A \subseteq \bigcup \calU$.
Suppose that for every $U \in \calU$, $\abs{A \cap U} \leq 1$,
and that $A$ dissectable.
Then $\calU$ is splittable.
\vspace{-2.3mm}
\end{prop}

\noindent
{\bf Proof } Trivial.\hfill\proofend

\begin{prop}\label{metr-bldr-p1.4}
Let $\rpair{X}{G}$ be a space-group pair,
and let $\,\calU$ denote a regular open family in $X$.
We define $\psi_{\fss{cnvrg}}(\calU)$
to be the following property of $\calU$.
\newline
\num{C1} $\calU$ is infinite and pairwise disjoint.
\newline
\num{C2} $\calU$ is absorbable.
\newline
\num{C3} $\calU$ is not splittable.

Let $\rpair{X}{G} \in K$ and $\calU \subseteq \fs{Ro}(X)$.

The following are equivalent.
\newline
\num{1} $\calU$ is pairwise disjoint and convergent,
and $\,\lim \calU \in \fs{D}(X,G)$.
\newline
\num{2} $\calU$ satisfies $\psi_{\fss{cnvrg}}$.

\vspace{-2.3mm}
\end{prop}

{\bf Proof }
$(1) \Rightarrow (2)$ \ Let $\calU$ be pairwise disjoint and convergent,
and\break
$\lim \calU \in D(X,G)$. Clearly, $\calU$ fulfills (C1).

Let $x = \lim \calU$.
There is
$V \in \fs{Ro}(X) \setminus \sngltn{\emptyset}$
such that $\setm{g(x)}{g \in G} \cap V$ is dense in $V$.
Let $W \subseteq V$ and $W \in \fs{Ro}(X) \setminus \sngltn{\emptyset}$.
Then there is $g \in G$ such that $g(x) \in W$.
So $\setm{U \in \calU}{g(U) \not\subseteq W}$ is finite.
Hence $\calU$ satisfies (C2).

Let $S \in \fs{Ro}(X) \setminus \sngltn{\emptyset}$.
We show that $S$ does not split $\calU$.
There is $W' \subseteq S$
such that $W' \in \fs{Ro}(X) \setminus \sngltn{\emptyset}$
and $\abs{\setm{U \in \calU}{U \cap W' \neq \emptyset}} \leq 1$.
Let $W \in \fs{Ro}(X) \setminus \sngltn{\emptyset}$ be such that
$\fs{cl}(W) \subseteq W'$.
Let $g \in G$, and suppose that
$\calV \eqdf \setm{V' \in \calU}{g(V') \cap W \neq \emptyset}$
is infinite.
Since $x = \lim \calV$, $g(x) = \lim g(\calV)$.
Hence
$g(x) \in \fs{cl}(W) \subseteq W'$.
So $\setm{V' \in \calU}{g(V') \not\subseteq W'}$ is finite.
Thus
\newline
$\setm{V' \in \calU}
{\mbox{there is } V'' \in \fs{Ro}(X) \setminus \sngltn{\emptyset}
\mbox{ such that } V'' \subseteq V' \mbox{ and }
g \restriction V'' = \fs{Id}}$\break
is finite.
Hence $S$ does not split $\calU$. So (C3) holds.

$(2) \Rightarrow (1)$ \
Suppose that $\kern2.8pt\calU$ satisfies $\psi_{\fss{cnvrg}}$.
So there is $V \in \fs{Ro}(X) \setminus \sngltn{\emptyset}$
such that $V$ absorbs $\calU$.
Let $S'$ be an open nonempty flexible set which intersects $V$ and
$S = S' \cap V$.
Then $\emptyset \neq S \subseteq V$ and $S$ is flexible.
Let\break
$g \in G$ be such that
$\setm{U \in \calU}{g(U) \not\subseteq S}$ is finite.
Let
$\calU' = \setm{g(U)}{U \in \calU \mbox{ and } g(U) \subseteq S}$.
Then $\calU'$ satisfies $\psi_{\fss{cnvrg}}$
and $\fs{acc}(\calU') = g(\fs{acc}(\calU))$.
So we may replace $\calU$ by $\calU'$. That is, we may assume that
for every $U \in \calU$, $U \subseteq S$.

We show that
$\emptyset \neq \fs{acc}(\calU) \subseteq D(X,G)$.
Let $D \subseteq S$ be a dense subset
such that every infinite discrete subset of $D$ is dissectable.
Suppose by contradiction that $\fs{acc}(\calU) = \emptyset$.
Let $A \subseteq D \cap \bigcup \calU$ be an infinite set
such that for every $U \in \calU$, $\abs{A \cap U} \leq 1$.
Then $\fs{acc}(A) = \emptyset$. Hence $A$ is dissectable.
By Proposition \ref{p1.3}, $\calU$ is splittabe. A contradiction.
So $\fs{acc}(\calU) \neq \emptyset$.

Let $x \in \fs{acc}(\calU)$.
Recall that $V$ absorbs $\calU$.
Let $W \subseteq V$ be a nonempty open set. Choose $W'$ such that
$W'$ is a nonempty regular open set and $\fs{cl}(W') \subseteq W$.
There is $g \in G$ such that
$\calA \eqdf \setm{U \in \calU}{g(U) \not\subseteq W'}$ is finite.
Since $x \in \fs{acc}(\calU \setminus \calA)$,
$g(x) \in \fs{acc}(g(\calU \setminus \calA)) \subseteq
\fs{cl}(W') \subseteq W$.
Hence $\setm{g(x)}{g \in G}$ is dense in $V$.
We have shown that $x \in D(X,G)$.

Let $x \in \fs{acc}(\calU)$. Suppose by contradiction that
$\calU$ does not converge to $x$.
Let $T \in \fs{Nbr}(x)$ be such that
$\setm{U \in \calU}{U \not\subseteq T}$ is infinite.
Let $A \subseteq (D \cap \bigcup \calU) \setminus T$
be an infinite set such that
for every $U \in \calU$, $\abs{A \cap U} \leq 1$.
Either (i) $A$ contains an infinite discrete set,
or (ii) $A$ contains the range of a $\onetoone$ convergent sequence.
Suppose that (i) happens. Since $A \subseteq D$, $A$ is dissectable.
By Proposition \ref{p1.3}, $\calU$ is splittable. A contradiction.

Suppose that (ii) happens. We may assume that $A$ is the range of
a $\onetoone$ convergent sequence.
Let $y = \lim A$.
Hence $y \in \fs{acc}(\calU)$ and by the previous argument,
$y \in D(X,G)$.
So either $\fs{DF}_{X,G}(x,y)$ holds or $\fs{DF}_{X,G}(y,x)$ holds.
We may assume that $\fs{DF}_{X,G}(x,y)$ holds.
So there is $V' \in \fs{Ro}(X) \setminus \sngltn{\emptyset}$ such that
$\setm{g(x)}{g \in G \mbox{ and for some } R \in \fs{Nbr}(y),
\ g \nrestriction R = \fs{Id}}$ is dense in $V'$.
We show that $V'$ splits $\calU$.
Let $W \subseteq V'$ be open and nonempty.
There are $g \in G$ and $R \in \fs{Nbr}(y)$ such that
$g(x) \in W$ and $g \nrestriction R = \fs{Id}$.
Since $x \in \fs{acc}(\calU)$ and $W \in \fs{Nbr}(g(x))$,
\hbox{$\setm{U \in \calU}{g(U) \cap W \neq \emptyset}$}
is infinite.\break
Since $A \cap R$ is infinite,
$\setm{U \in \calU}{U \cap R \neq \emptyset}$ is infinite.
Hence
\newline
$\setm{U \in \calU}
{\mbox{there is } U' \in \fs{Ro}(X) \setminus \sngltn{\emptyset}
\mbox{ such that } U' \subseteq U
\mbox{ and } g \restriction U' = \fs{Id}}$\break
is infinite. That is, $V'$ splits $\calU$. A contradiction.
So $\lim \calU = x$.
\rule{0pt}{0pt}\hfill\proofend

\begin{prop}\label{metr-bldr-p1.6}
Let $\psi^1_{\fss{eq}}(\calU,\calV)$ be the following property.
\newline
There does not exist a nonempty regular open set $S$
such that:
For every $T \in \fs{Ro}(X) \setminus \sngltn{\emptyset}$:
if $T \subseteq S$, then there is $g \in G$
such that $\setm{U \in \calU}{g(U) \not\subseteq T}$ is finite,
and $\setm{V \in \calV}{g \restriction V \neq \fs{Id}}$ is finite.
Let
\newline\centerline{
$\psi_{\fss{eq}}(\calU,\calV) \equiv
\psi^1_{\fss{eq}}(\calU,\calV) \wedge \psi_{\fss{eq}}^1(\calV,\calU)$.
}

Let $\calU,\calV$ be convergent pairwise disjoint regular open
families such that $\lim \calU, \lim \calV \in D(X,G)$.
Then $\psi_{\fss{eq}}(\calU,\calV)$ holds
iff $\,\lim \calU = \lim \calV$.
\vspace{-2.3mm}
\end{prop}

\noindent
{\bf Proof }
Suppose that $\lim \calU = \lim \calV$.
Let $S \in \fs{Ro}(X) \setminus \sngltn{\emptyset}$.
Let $x = \lim \calU$.
Choose $T \in \fs{Ro}(X) \setminus \sngltn{\emptyset}$ such that
$T \subseteq S$ and $x \not\in \fs{cl}(T)$.
Let $g \in G$\break
be such that $\calA \eqdf \setm{U \in \calU}{g(U) \not\subseteq T}$
is finite.
Since $x = \lim (\calU \setminus \calA)$, $g(x) \in \fs{cl}(T)$.
So $g(x) \neq x$.
Let $Q,R$ be pairwise disjoint neighborhoods of $x$ and $g(x)$
respectively.
Since $x = \lim \calV$, $g(x) = \lim g(\calV)$.
Hence
\newline
{\thickmuskip=2mu \medmuskip=1mu \thinmuskip=1mu
$\setm{V \in \calV}{g(V) \not\subseteq R}$ is finite.
Since $\lim \calV = x$ and $Q \in \fs{Nbr}(x)$,
\hbox{$\setm{V \in \calV}{V \not\subseteq Q}$}
}
is finite.
Hence $g(V) \cap V = \emptyset$ for all but finitely many members of
$\calV$.
So $\setm{V \in \calV}{g \nrestriction V \neq \fs{Id}}$ is infinite.
Hence $\rpair{\calU}{\calV}$ fulfills $\psi^1_{\fss{eq}}$.
Similarly, $\rpair{\calV}{\calU}$ fulfills $\psi^1_{\fss{eq}}$.
So $\rpair{\calU}{\calV}$ fulfills $\psi_{\fss{eq}}$.

Let $\calU,\calV$ be convergent pairwise disjoint regular open
families such that $\lim \calU, \lim \calV \in D(X,G)$.
Suppose that $x = \lim \calU \neq \lim \calV = y$.
By (Q2) of Definition \ref{d1.1}(f),
$\fs{DF}_{X,G}(x,y)$ holds or $\fs{DF}_{X,G}(y,x)$ holds.
We may assume that $\fs{DF}_{X,G}(x,y)$ holds.
Let $S$ be a regular open set such that
\newline
$\setm{g(x)}{g \in G \mbox{ and for some } R \in \fs{Nbr}(y) \mbox{, }
g \nrestriction R = \fs{Id}}$ is dense in $S$.
Let $T \subseteq S$ be regular open and nonempty.
There is $g \in G$ such that $g(x) \in T$
and for some $R \in \fs{Nbr}(y)$, $g \nrestriction R = \fs{Id}$.
Then $\setm{U \in \calU}{g(U) \not\subseteq T}$ is finite.
Let $\calV' = \setm{V \in \calV}{V \subseteq R}$.
Then $\calV \setminus \calV'$ is finite,
and for every $V \in \calV'$, $g \nrestriction V = \fs{Id}$.
That is, $\setm{V \in \calV}{g \nrestriction V \neq \fs{Id}}$
is finite.
Hence $\rpair{\calU}{\calV}$ does not fulfill $\psi_{\fss{eq}}$.
\hfill\proofend

\begin{prop}\label{metr-bldr-p5.31}
Let $\psi_{\in}(\calU,V)$ be the following property.
\newline
For every $\calW$: if $\psi_{\fss{cnvrg}}(\calW)$ and
$\psi_{\fss{eq}}(\calU,\calW)$,
then $\setm{W \in \calW}{W \not\subseteq V}$ is finite.
\newline
Then for every $\calU$ satisfying $\psi_{\fss{cnvrg}}$
and $V \in \fs{Ro}(X)$:
$\psi_{\in}(\calU,V)$ holds iff $\lim \calU \in V$.
\vspace{-2.3mm}
\end{prop}

{\bf Proof } It is obvious that if $\lim \calU \in V$,
then $\psi_{\in}(\calU,V)$ holds.
Suppose that $\lim \calU \not\in V$. Let $x = \lim \calU$.
Since $V \in \fs{Ro}(X)$, $x \not\in \fs{int}(\fs{cl}(V)))$.
So $x \in \fs{cl}(X \setminus \fs{cl}(V))$. Hence there is
a regular open pairwise disjoint family $\calW$
such that for every $W \in \calW$, $W \subseteq X \setminus \fs{cl}(V)$
and $\lim \calW = x$.
So $\psi_{\fss{eq}}(\calU,\calW)$ holds and
$\setm{W \in \calW}{W \not\subseteq V}$ is not finite.
So $\psi_{\in}(\calU,V)$ does not hold.\break
\rule{10cm}{0pt}
\hfill\proofend

\noindent
{\bf Proof of Theorem \ref{t1.2} }
Let $\rpair{X}{G},\rpair{Y}{H} \in K$ and $\iso{\varphi}{G}{H}$.
By Theorem \ref{t2.3}, there is
$\iso{\eta}{\bfs{Ro}(X)}{\bfs{Ro}(X)}$ such that
$\varphi(g) = \eta \scirc g \scirc \eta\inverse$ for every $g \in G$.
That is, for every $g \in G$ and $U \in \fs{Ro}(X)$,
$\varphi(g)(\eta(U)) = \eta(g(U))$.

We define $\fnn{\tau}{D(X,G)}{D(Y,H)}$.
Let $x \in D(X,G)$.
Let $\calU$ be a pairwise disjoint regular open family such that
$\lim \calU = x$.
By Proposition~\ref{metr-bldr-p1.4}, $\calU$ satisfies
$\psi_{\fss{cnvrg}}$.
So $\eta(\calU)$ satisfies $\psi_{\fss{cnvrg}}$.
By Proposition~\ref{metr-bldr-p1.4}, $\eta(\calU)$ converges to
a member of $D(Y,H)$.
Define $\tau(x) = \lim \eta(\calU)$.
If $\calV$ is another pairwise disjoint regular open family
converging to $x$, then by Proposition~\ref{metr-bldr-p1.6},
$\psi_{\fss{eq}}(\calU,\calV)$ holds.
So $\psi_{\fss{eq}}(\eta(\calU),\eta(\calV))$ holds.
Hence by Proposition~\ref{metr-bldr-p1.6},
$\lim \eta(\calU) = \lim \eta(\calV)$.
So the definition of $\tau$ is independent of the choice of $\calU$.
It follows easily that $\tau$ is a bijection between
$D(X,G)$ and $D(Y,H)$.

We show that $\iso{\tau}{D(X,G)}{D(Y,H)}$.
It suffices to show that for every $x \in D(X,G)$
and $V \in \fs{Ro}(X)$,
$x \in V$ iff $\tau(x) \in \eta(V)$.
Let $\calU$ be a pairwise disjoint regular open family
converging to $x$. By Proposition~\ref{metr-bldr-p5.31},\break
$\psi_{\in}(\calU,V)$ holds iff $x \in V$.
Also,
$\psi_{\in}(\calU,V)$ holds iff $\psi_{\in}(\eta(\calU),\eta(V))$
holds iff $\lim \eta(\calU) \in \eta(V)$.
So $x \in V$ iff $\lim \eta(\calU) \in \eta(V)$.
That is, $x \in V$ iff $\tau(x) \in \eta(V)$.

It is left to the reader to check that for every $g \in G$,
\newline\centerline{
\kern3cm
$\varphi(g) \nrestriction D(Y,H) =
\tau \scirc (g \nrestriction D(X,G)) \scirc \tau\inverse$.
\hfill\proofend
}

\section{The non-shrinking property of the semi-groups $\itGamma_i(E)$}
\label{s3}

In this section we consider two subsemigroups of
$\rpair{\setm{f}{\fnn{f}{X}{X}}}{\scirc}$,
where $X$ is an infinite dimensional normed space:
$\itGamma_1(X)$ and $\itGamma_2(X)$.
By their definition, $\itGamma_1(X) \subseteq \itGamma_2(X)$.

We prove that (1) for every infinite dimensional normed space $X$,
$\itGamma_1(X)$ is non-shrinking.
We also show that (2) if $X$ is an infinite dimensional normed space,
and $X$ has an infinite-dimensional
subspace $Y$, such that the completion of $Y$ does not have a subspace
isomorphic to the space $c_0$ (of real sequences converging to $0$),
then $\itGamma_2(X)$ is non-shrinking.

Finally we prove that (3) $\itGamma_2(c_0)$ is not non-shrinking.

The proofs of (1) and (2) have a very similar structure.
In both proofs we need to show that the unit sphere of $X$
does not have a cover with certain properties.
But for spaces which do not embed $c_0$ a stronger statement of this
kind can be proved.
The proof of (2) relies on the fact (Lemma \ref{l3.8}) that
if $X$ is an infinite dimensional Banach space,
and does not have a subspace isomorphic to $c_0$
then $X$ has the property:
\newline
($*$) The sphere of $X$ does not have a locally finite cover
consisting of weakly closed sets which do not contain $0$.

This fact is not true for a general infinite dimensional Banach space,
and in fact, $c_0$ does not have this property.

The following fact (Corollary \ref{p3.6})
replaces ($*$) when proving (1).
Every infinite dimensional normed space $X$ has the property:
\newline
($*$$*$) The sphere of $X$ does not have a locally finite
cover with finite order consisting of closed convex sets
which do not contain $0$.
\newline
A set $\calS$ of subsets of a set $X$ has finite order $k \in \bbN$,
if
\newline
$\max_{x \in X} \abs{\setm{S \in \calS}{x \in S}} = k$.

Let $E$ be a normed space.
Denote by $\overE$ the completion of $E$.
We define a set of functions $\whatGamma_2 = \whatGamma_2(E)$
from $E$ to $E$.
A function $\fnn{h}{E}{E}$
belongs to $\whatGamma_2$ if there is a set of pairs
\newline
$\setm{\rpair{\whatS_j}{E_j}}{j \in J}$ such that:
\newline
(1) Every $\whatS_j$ is a weakly closed subset of $\overE$.
\newline
(2) $\setm{\whatS_j)}{j \in J}$ is locally finite
in $\overE$.
\newline
(3) $\bigcup \setm{\whatS_j}{j \in J}$ is bounded.
\newline
(4) For every $j \in J$ and $x \in \whatS_j \cap E$,
\newline\centerline{
$h(x) - x \in E_j$,
}
and for every
$x \in E \setminus \bigcup \setm{\whatS_j}{j \in J}$,
\newline\centerline{
$h(x) = x$.
}
For every $h \in \whatGamma_2$ we pick a family
$\setm{\rpair{\whatS_j}{E_j}}{j \in J}$ satisfying the above
and in which $\whatS_i \neq \whatS_j$ for every distinct $i,j \in J$.
We then denote $J$, $\whatS_j$, $E_j$ by\break
$J^h$, $\whatS_J^h$ and $E_j^h$ respectively.
If $X$ is a metric space, $x \in X$ and $r > 0$,
define $B_X(x,r) = \setm{y \in X}{d(x,y) \leq r}$.
For a normed space $E$ denote $B_E(0,1)$ by $B_E$.

Let $\itGamma_2 = \itGamma_2(E)$ be the semigroup generated by
$\whatGamma_2$,
that is, the closure of $\whatGamma_2$ under composition.

Note that $\whatGamma_2$ contains discontinuous functions.

Let $\calS$ be a set of subsets of a metric space $\rpair{X}{d}$
and $r > 0$.
$\calS$ is\break
{\it $r$-separated} if for every distinct $S,T \in \calS$,
$d(S,T) > r$.
We say that $\calS$ is {\it separated} if for some $r > 0$,
$\calS$ is $r$-separated.

We define $\whatGamma_1 = \whatGamma_1(E)$.
A function $\fnn{h}{E}{E}$
belongs to $\whatGamma_1$ if there is a set of pairs
$\setm{\rpair{S_j}{E_j}}{j \in J}$ such that:
\newline
(1) Every $S_j$ is a closed convex subset of $E$.
\newline
(2) $\setm{S_j}{j \in J}$ is separated.
\newline
(3) $\bigcup \setm{S_j}{j \in J}$ is bounded.
\newline
(4) For every $j \in J$ and $x \in S_j$,
\newline\centerline{
$h(x) - x \in E_j$,
}
and for every
$x \in E \setminus \bigcup \setm{S_j}{j \in J}$,
\newline\centerline{
$h(x) = x$.
}
Let $\itGamma_1 = \itGamma_1(E)$ be the semigroup generated by
$\whatGamma_1$.

Note that $\whatGamma_1 \subseteq \whatGamma_2$.
For $h  \in \whatG_1$,
define $\whatS^h_j = \fs{cl}_{\oversE}(S^h_j)$.
Then
\newline
$\setm{\rpair{\whatS^h_j}{E^h_j}}{j \in J^h}$ fulfills
Clauses (1)\,-\,(4) in the definition of $\itGamma_2$.

The theorem below is the final goal of this section.

\begin{theorem}\label{t3.1}
\num{a} For every infinite dimensional normed space $E$,
$\itGamma_1(E)$ is non-shrinking.

\num{b} Let $E$ be an infinite dimensional normed space.
Suppose that $E$ contains an infinite dimensional subspace $F$
such that $c_0$ is not isomorphic to a subspace of the completion
of $F$. Then $\itGamma_2(E)$ is non-shrinking.
\end{theorem}

For a function $g$, $\fs{supp}(g)$ is defined as
$\setm{x \in \fs{Dom}(g)}{g(x) \neq x}$.
Clearly,
$\fs{supp}(h_n \scirc \ldots \scirc h_1) \subseteq
\bigcup_{i = 1}^n \fs{supp}(h_i)$.
So for every $g \in \itGamma_2$, $\fs{supp}(g)$ is bounded.

The definition of the non-shrinking property relies on the choice
of the metric on $X$.
We observe that the non-shrinkingness of $\itGamma_1$ and
$\itGamma_2$ persists when one norm on $E$ is replaced by
an equivalent norm.

\begin{prop}\label{p3.2}
\num{a} If $\itGamma_2$ is not non-shrinking,
then for every $s \in (0,1)$
there are $g \in \itGamma_2$ and a finite
subset $\sigma \subseteq E$ such that
\newline
$g(B_E) \subseteq \bigcup_{x \in \sigma} (x + s \ncdot B_E)$.

\num{b} The same holds for $\itGamma_1$.
That is, if $\itGamma_1$ is not non-shrinking,
then for every $s \in (0,1)$ there are
$g \in \itGamma_1$ and a finite
subset $\sigma \subseteq E$ such that
$g(B_E) \subseteq \bigcup_{x \in \sigma} (x + s \ncdot B_E)$.

\num{c} Let $(E,\norm{.})$ be a normed space,
and $\3norm{.}$ be an equivalent norm on $E$.
If $\itGamma_2$ is non-shrinking with respect to $\norm{.}$,
then it is non-shrinking with respect to $\3norm{.}$.

\num{d} The same holds for $\itGamma_1$.
\vspace{-2.3mm}
\end{prop}

\noindent
{\bf Proof }
(a) Suppose that $g \in \itGamma_2$, $r \in (0,1)$
and $\sigma \subseteq E$ are such that
$\sigma$ is finite and
$g(B_E)  \subseteq \bigcup_{x \in \sigma} (x + r \ncdot B_E)$.

For $v \in E$ and $A \subseteq E$ let
$\bfs{tr}_v$ be the function $x \mapsto x + v$
and
\newline
$\bfs{tr}_{v,A} =
\bfs{tr}_v \nrestriction A \cup \fs{Id} \nrestriction (E \setminus A)$.
So if $A$ is bounded, then $\bfs{tr}_{v,A} \in \itGamma_1$.

For $v \in E$ and $\lambda > 0$ let
$g_{\lambda,v}$ be the function acting on
$v + \lambda \ncdot B_E$ in the same way that $g$ acts
on $B_E$.
That is,
$g_{\lambda,v} =
\bfs{af}_{\lambda,v} \scirc g \scirc (\bfs{af}_{\lambda,v})\inverse$,
where $\bfs{af}_{\lambda,v}$ is the affine function:
$x \mapsto \lambda x + v$.
Clearly, $g_{\lambda,v} \in \itGamma_2$ and
\newline
$g_{\lambda,v}(v + \lambda \ncdot B_E) \subseteq
\bigcup_{x \in \sigma} (x + v + r \lambda \ncdot B_E)$.

Let $\sigma = \fsetn{x_1}{x_k}$.
Choose $\setm{v_i}{i \leq k}$
such that for every $i < j \leq k$,
\newline
(1)
$\fs{supp}(g_{v_i,r}) \cap (v_j + r \ncdot B_E) = \emptyset$;
\newline
(2)
$g_{v_i,r}(v_i + r \ncdot B_E) \cap \fs{supp}(g_{v_j,r}) =
\emptyset$;
\newline
(3)
$v_i + r \ncdot B_E \cap
\bigcup_{\ell = 1}^k (x_{\ell} + r \ncdot B_E) =
\emptyset$.
\newline
Let
$A_i = x_i + r \ncdot B_E \setminus
\bigcup_{\ell < i} (x_{\ell} + r \ncdot B_E)$,
$f_i = \bfs{tr}_{v_i - x_i,A_i}$
and
$f = f_k \scirc \ldots \scirc f_1$.
Relying on (3),
\newline
(3.1)
\ For every $i \leq k$, \indent
$f(x_i + r \ncdot B_E \setminus
\bigcup_{\ell < i} (x_{\ell} + r \ncdot B_E))
\subseteq
v_i + r \ncdot B_E$.
\newline
Let
\newline\centerline{
$h = g_{v_k,r} \scirc \ldots \scirc g_{v_1,r}$
\ and \ $\hath = h \scirc f$.
}
We check that
\newline\centerline{
$\hath \scirc g(B_E) \subseteq
\bigcup \setm{x_i + v_j + r^2 \ncdot B}{1 \leq i \leq k,\ %
1 \leq j \leq k}$.
}
For every $i \leq k$,
$g_{v_i,r}(v_i + r \ncdot B_E) \subseteq
\bigcup_{j \leq k} (x_j + v_i + r^2 \ncdot B_E)$.
For every $\ell < i$,
$g_{v_{\ell},r} \restriction (v_i + r \ncdot B_E) = \fs{Id}$.
This follows from (1).
So
\newline
(3.2) \
$g_{v_i,r} \scirc \ldots \scirc g_{v_1,r}
(v_i + r \ncdot B_E) =
g_{v_i,r}(v_i + r \ncdot B_E) \subseteq
\bigcup_{j \leq k} (x_j + v_i + r^2 \ncdot B_E)$.
\newline
Clause (2) together with the equality in (3.2) imply that
\newline
(3.3) For every $i \leq k$
\newline\centerline{
$h(v_i + r \ncdot B_E) =
g_{v_i,r} \scirc \ldots \scirc g_{v_1,r}(v_i + r \ncdot B_E)
\subseteq
\bigcup_{j \leq k} (x_j + v_i + r^2 \ncdot B_E)$.
}
By (3.1) and (3.3),
\newline
(3.4) \ For every $i \leq k$, \indent
$\hath(A_i) \subseteq
\bigcup_{j \leq k} (x_j + v_i + r^2 \ncdot B_E)$.
\newline
Since
$\bigcup_{i \leq k} A_i =
\bigcup_{i \leq k}(x_i + r \ncdot B_E)$,
\newline
(3.5)
\centerline{
$\hath(\bigcup_{i \leq k}(x_i + r \ncdot B_E)) \subseteq
\bigcup_{j \leq k} (x_j + v_i + r^2 \ncdot B_E)$.
}
It follows that
$\hath \scirc g(B_E) \subseteq
\bigcup_{j \leq k} (x_j + v_i + r^2 \ncdot B_E)$.
Hence it is contained in the union of finitely many balls with radius
$r^2$.

Let $s \in (0,1)$. Iterating the above construction sufficiently many
times, one obtains $g_s \in \itGamma_2$ such that $g_s(B_E)$
is contained in the union of finitely many balls with radius $\leq s$.

(b) Note that if in the above construction $g \in \itGamma_1$,
then $\hath \in \itGamma_1$.

(c) Let $B'$ be the unit ball of $(E,\3norm{.})$.
Let $K > 1$ be such that $B' \subseteq K \ncdot B_E$
and $B_E \subseteq K \ncdot B'$.
Let $S > K^2$.
By Part (a), there is $g \in \itGamma_2$ such that
$g(K \ncdot B_E)$ is contained in a finite union
$\bigcup_{i \leq k} D_i$ of $\norm{.}$-balls with radius $\frac{K}{S}$.
So
$g(K \ncdot B_E) \subseteq \bigcup_{i \leq k} D_i'$,
where each $D_i'$ is a $\3norm{.}$-ball with radius $\frac{K^2}{S}$.
So $g(B') \subseteq \bigcup_{i \leq k} D_i'$,
and the $D_i'$\raise1pt\hbox{-}s are $\3norm{.}$-balls with radius
$< 1$.

(d) Part (d) follows from Part (b) in the same way that
Part (c) follows from Part (a).
\hfill\proofend

Suppose that $L,M$ are linear subspaces of a vector space $X$
and
\newline
$L \cap M = \sngltn{0}$.
The function $\fnn{P}{M + L}{L}$ defined by
\newline\centerline{
$P(v + u) = u$, \ \,$v \in M$, $u \in L$
}
is called the projection of $\rpair{L}{M}$.

\begin{lemma}\label{l3.3}
Let $X$ be a separable normed space.
Then there is an equivalent norm $\3norm{.}$ on $X$ such that

\noindent
\num{*} for any finite-dimensional subspace
$M \subseteq X$ and any $\varepsilon > 0$
there is a finite-codimensional subspace $L \subseteq X$
such that $M \cap L = \sngltn{0}$ and if
\newline
$\fnn{P}{M + L}{L}$
is the projection of $\rpair{L}{M}$, then $\3norm{P} < 1 + \varepsilon$.
\end{lemma}

\noindent
{\bf Proof.}
Let $Y$ be any separable Banach space with a basis which contains
$X$ isomorphically. For instance, take $Y = C([0,1])$.
Let $\setm{e_i}{i \in \bbN}$ be a basis for $Y$ such that for every
$i \in \bbN$, $\norm{e_i} = 1$.
Let $\setm{e_i^*}{i \in \bbN} \subseteq Y^*$
be the biorthogonal sequence for $\setm{e_i}{i \in \bbN}$.
For $n = 0,1,\ldots$ define the operators from $Y$ to $Y$
$$
S_n x = \sum_{i = 1}^n e_i^*(x)e_i \mbox{\ \ and\ \ }
R_n x = x - S_n x,
$$
and put
$$\3norm{x}=
\sup\setm{||S_n x||,\; ||R_n x||\;}{\; n = 0,1,...} .$$
So $S_0$ is the $0$-operator and $R_0 = \fs{Id}$.
The norm $\3norm{.}$ has the following properties.

(P1) \ $\3norm{.}$ is equivalent to the original norm on $Y$.

(P2) \ $\3norm{S_n} = \3norm{R_n} = 1,\; n = 1,2,\ldots$

{\bf Proof of P1 }
For every $x \in X$,
$\3norm{x} \geq \3norm{R_0 x} = \norm{x}$.

The fact that $\setm{x_n}{n \in \bbN}$ is a basic sequence is
equivalent to the existence of $C > 0$
such that for every $x \in X$ and $n \in \bbN$,
$\norm{S_n x} \leq C \norm{x}$.
Hence for every $n \in \bbN$,
$\norm{R_n x} = \norm{x - S_n x} \leq (1 + C) \ncdot \norm{x}$.
Since $\3norm{.}$ is the supremum of the above,
$\3norm{x} \leq (1 + C) \ncdot \norm{x}$.

{\bf Proof of P2 }
Since for $n = 1,2,\ldots, \ S_n$ and $R_n$ are nonzero projections,
their norm must be $\geq 1$.
It follows trivially from the definition of $\3norm{.}$
that $\3norm{R_n},\3norm{S_n} \leq 1$.

We claim that the norm
$\3norm{.}$ restricted to $X$ has Property ($*$).
\newline
Let $Y_m = [e_i]_{i = 1}^m$ and
$Y^m = [e_i]_{i = m + 1}^{\infty}$.
Then $\fnn{S_m}{Y^m + Y_m}{Y_m}$ is the projection of
$\rpair{Y_m}{Y^m}$ and
$\fnn{R_m}{Y_m + Y^m}{Y^m}$ is the projection of
$\rpair{Y^m}{Y_m}$.
Since $\3norm{S_m} = \3norm{R_m} = 1$,
for every $u \in Y_m$ and $v \in Y^m$,
$\3norm{u},\3norm{v} \leq \3norm{u + v}$.

Let $M \subseteq X$ be a finite dimensional subspace of $X$
and $\varepsilon > 0$.
We may assume that $\varepsilon < 1$.
Let $m \in \bbN$ be so large that for every $x \in S_M$,
$\inf(\setm{\3norm{x - y}}{y \in Y_m}) \leq \frac{\varepsilon}{3}$.
Such an $m$ exists, since $S_M$ is compact and
$\bigcup_{i = 1}^{\infty} Y_m$ is dense in Y.
We check that $M \cap Y^m = \sngltn{0}$.
If this is not so, there is $x \in S_M \cap Y^m$.
Let $y \in Y_m$ be such that
$\3norm{x - y} \leq \frac{\varepsilon}{3}$.
Hence $\3norm{y} \geq 1 - \frac{\varepsilon}{3}$.
So for every $z \in Y^m$,
$\3norm{y - z} \geq  1 - \frac{\varepsilon}{3}$.
In particular, $\3norm{y - x} \geq  1 - \frac{\varepsilon}{3}$.
This is impossible since $\frac{\varepsilon}{3} \leq \third$.

Let $\fnn{Q}{M + Y^m}{Y^m}$ be the projection of
$\rpair{Y^m}{M}$.
We show that $\3norm{Q} \leq 1 + \varepsilon$.
Let $x \in M + Y^m$ and $\3norm{x} = 1$.
So for some $u \in M$ and $v \in Y^m$,
$x = u + v$.
Let $w \in Y^m$ be such that
$\3norm{u - w} \leq \frac{\varepsilon}{3} \ncdot \3norm{u}$.
So
\begin{equation}\label{m1}
\3norm{v} \leq \3norm{v + w} \leq \3norm{v + u} + \3norm{u - w} \leq
\3norm{v + u} + \frac{\varepsilon}{3} \cdot \3norm{u}.
\end{equation}
\begin{equation}\label{m2}
\3norm{u} \leq \3norm{u + v} + \3norm{-v} =
1 + \3norm{v} \leq 1 + \3norm{v + w} \leq
\3norm{v + u} + \frac{\varepsilon}{3} \ncdot \norm{u}.
\end{equation}
Hence
\begin{equation}\label{m3}
(1 - \frac{\varepsilon}{3}) \ncdot \3norm{u} \leq \3norm{v + u}.
\end{equation}
That is,
\begin{equation}\label{m4}
\3norm{u} \leq \frac{3}{3 - \varepsilon} \ncdot \3norm{v + u}.
\end{equation}
By (\ref{m1}) and (\ref{m4}),
\begin{equation}\label{m5}
\3norm{v} \leq \3norm{v + u} + \frac{\varepsilon}{3} \ncdot
\frac{3}{3 - \varepsilon} \ncdot \3norm{v + u} =
(1 + \frac{\varepsilon}{3 - \varepsilon}) \ncdot \3norm{v + u} \leq
(1 + \varepsilon) \ncdot \3norm{v + u}.
\end{equation}

So $\3norm{Q} \leq 1 + \varepsilon$.
Let $L = Y^m \cap X$.
Since $Y^m$ has finite codimension in~$Y$,
$L$ has finite codimension in $X$.
Clearly, $M \cap L = \sngltn{0}$.
The projection $P$ of $\rpair{L}{M}$ is the restriction of $Q$
to $M + L$.
So $\3norm{P} \leq 1 + \varepsilon$.\hfill\proofend

Let $\calS$ be a set of sets and $k \in \bbN$.
$\calS$ has {\it finite order}, if there is $k \in \bbN$ such that
for every $x$, $\abs{\setm{S \in \calS}{x \in S}} \leq k$.
The {\it order} of $\calS$ is the minimal such $k$.

\begin{lemma}\label{l3.4}
Let $M \subseteq X$ be a finite dimensional
subspace of a normed space $X$ and $L \subseteq X$ be a closed
subspace which is a complement of $M$ in $X$.
Let $\fnn{P}{X}{L}$ be the projection of $\rpair{L}{M}$.

\num{a} Let $A \subseteq X$ be a bounded closed set.
Then $P(A)$ is closed.

\num{b} Let $\calS$
be a locally finite family of subsets of $X$ such that
$\bigcup_{S \in \calS} S$ is bounded.
Then $\setm{P(S)}{S \in \calS}$ is a locally finite family.

\num{c} Let $\calS$ be a separated family of subsets of $X$ such that
$\bigcup_{S \in \calS} S$ is bounded.
Then $\setm{P(S)}{S \in \calS}$ has finite order.
\end{lemma}

\noindent
{\bf Proof }
(a) Let $x = \lim x_n ,\; x_n \in P(A)$.
Then for every $n$ there is $z_n \in A$ such that $Pz_n = x_n$.
Clearly, $z_n = x_n + y_n$, where $y_n \in M$.
Since $A$ is bounded,
$\setm{z_n}{n \in \bbN}$ is a bounded sequence.
$x_n$ is convergent and thus it is bounded. So $y_n = z_n - x_n$
is a bounded sequence.
Also, $\setm{y_n}{n \in \bbN} \subseteq M$,
and $M$ is finite-dimensional.
So $\setm{y_n}{n \in \bbN}$ has a convergent subsequence
$\setm{y_{n_k}}{k \in \bbN}$.
Denote $y  = \lim_k y_{n_k}$.

Then we have $\lim z_{n_k} = \lim_k \,(x_{n_k} + y_{n_k}) = x + y$.
Since $A$ is closed and $z_{n_k}\in A$,
it follows that $x + y \in A$. Clearly, $P(x + y) = x \in P(A)$.
This proves (a).

(b) Let $\calS$ be as in Part (b), and assume by contradiction that
$\setm{P(S)}{S \in \calS}$ is not locally finite.
Let $x$ be an accumulation point of $P(\calS)$.
That is, there is a $\onetoone$ sequence
$\setm{S_i}{i \in \bbN} \subseteq \calS$,
$x_i \in P(S_i)$ such that $\lim_i x_i = x$.
Let $y_i \in M$ be such that $y_i + x_i \in S_i$.
Since $\bigcup_{S \in \calS} S$ is bounded and the projection of
$\rpair{M}{L}$ is bounded,
$\setm{y_i}{i \in \bbN}$ is bounded.
Since $M$ is finite dimensional $\setm{y_i}{i \in \bbN}$ has a
convergent subsequence. Denote it by $\setm{y_{i_k}}{k \in \bbN}$.
Both $\sngltn{y_{i_k}}$ and $\sngltn{x_{i_k}}$ are convergent.
So $\sngltn{x_{i_k} + y_{i_k}}$ is convergent.
But $x_{i_k} + y_{i_k} \in S_{i_k}$, and the $S_{i_k}$'s are distinct.
So $\calS$ is not locally finite. A contradiction.
This proves Part (b).

(c) A subset $A$ of a metric space is $\varepsilon$-separated,
if the distance between every two distinct points of $A$ is
$\geq \varepsilon$.
There is an integer
$\ell = \ell(m,R,\varepsilon)$
such that for every $m$-dimensional normed space $E$
and an $\varepsilon$-separated subset $A \subseteq B_E(0,R)$,
$\abs{A} \leq \ell$.
For spaces $E$ isometric to $\bbR^m_{\infty}$ the number\break
$\ell' = \left(\left[\frac{R}{\varepsilon}\right] + 1 \right)^m$
fulfills the requirement.

Let $E$ be an $m$-dimensional normed space.
There is a norm $\3norm{.}$ on $E$ such that
$\rpair{E}{\3norm{.}}$ is isometric to $\bbR^m_{\infty}$
and for every $x \in E$,
$\frac{\3norm{x}}{m} \leq \norm{x} \leq m \3norm{x}$.

Let $A \subseteq B_E(0,R)$ be $\varepsilon$-separated.
Then $A \subseteq B_E^{\3norm{.}}(0,mR)$ and $A$ is\break
\rule{0pt}{18pt}$\frac{\varepsilon}{m}$-separated with respect to
$\3norm{.}$.
So
$\abs{A} \leq
\left(\left[\frac{m^2 R}{\varepsilon}\right] + 1 \right)^m$.

Suppose that $\fs{dim}(M) = m$ and let $\calS$ be as in Part (c).
Let $\delta > 0 $ be such that $\calS$ is $\delta$-separated.
Suppose that $\bigcup_{S \in \calS} S \subseteq B_X(0,R)$.
Let $Q$ be the projection of $\rpair{M}{L}$.
So
$Q(\bigcup_{S \in \calS} S) \subseteq B_M(0,\norm{Q} \ncdot R)$.
\hbox{Denote $R' = \norm{Q} \ncdot R$.}
Let $k \in \bbN$ be such that there are distinct
$S_1,\ldots,S_k \in \calS$ and $x$
with $ x \in \bigcap_{i = 1}^k P(S_i)$.
So there are $y_1,\ldots,y_k \in M$ such that
for every $i = 1,\ldots,k$, \ $y_i + x \in S_i$.
It follows that $\fsetn{y_1}{y_k} \subseteq B_M(0,R')$.
Also, for every $i < j \leq k$,
$\norm{y_i - y_j} = \norm{(x + y_i) - (x + y_j)} \geq \delta$.
So $k \leq \ell(m,R',\delta)$.
It follows that the order of $\calS$ is $\leq \ell(m,R',\delta)$.
\hfill\proofend

For a normed space $E$, $S_E$ denotes the unit sphere of $E$.

The next lemma is due to Michael Levin.

\begin{lemma} \label{l3.5} {\rm (M. Levin)}
Let $X$ be an $n$-dimensional normed space.
Suppose that $\calF$ is a finite set of closed subsets of $X$ such that
\newline
\num{1} $S_X = \bigcup_{F \in \calF} F$.
\newline
\num{2} For every $F \in \calF$ and $x \in F$, $-x \not\in F$.
\newline
Then the order of $\calF$ is $\geq \frac{n}{2} + 1$.
\vspace{-2.3mm}
\end{lemma}

\noindent
{\bf Proof }
Let $k$ be the order of $\calF$ and assume by contradiction that\break
$k < \frac{n}{2} + 1$.
There is $\varepsilon > 0$ such that the order of
$\calU \eqdf \setm{B_X(F,\varepsilon)}{F \in \calF}$ is $k$,
and for every $x \in U \in \calU$, \,$-x \not\in U$.
Let $\setm{\psi_U}{U \in \calU}$ be a partition of unity for $\calU$.
That is,
\newline
(1) For every $U \in \calU$,
$\fs{Dom}(\psi_U) = \bigcup_{V \in \calU} V$.
Denote $W = \bigcup_{V \in \calU} V$.
\newline
(2) For every $U \in \calU$ and $x \in U$, \ $\psi_U(x) > 0$.
\newline
(3) For every $U \in \calU$
and $x \in W \setminus U$, \ $\psi_U(x) = 0$.
\newline
(4) For every $x \in W$, \ $\sum_{U \in \calU} \psi_U(x) = 1$.

Let $K$ be the simplex whose vertices are the members of $\calU$,
and define $\fnn{f}{S_X}{K}$.
\newline\centerline{
$f(x) = (\psi_U(x))_{U \in \calU}$, \ $x \in S_E$.
}
Let $L = \fs{Rng}(f)$, and denote by $m$ the topological dimension of
$L$. If $z \in L$, then the number of coordinates of $z$
which are different from $0$ is $\leq k$.
That is, $L$ is contained in the
{\thickmuskip=2mu \medmuskip=1mu \thinmuskip=1mu
$(k - 1)$-dimensional skeleton of $K$.
}
So $m \leq k - 1$.
\newline
$2(k - 1) + 1 < 2(\frac{n}{2} + 1 - 1) + 1 = n + 1$.
That is, $2m + 1 \leq n$.
By \cite{E} Theorem 1.11.4 p.95,
there is an embedding $\fnn{g}{L}{\bbR^n}$.
It is obvious that if $x \in S_X$,
then $g \scirc f(x) \neq g \scirc f(-x)$.
This contradicts the theorem of Borsuk and Ulam which says that
for every continuous function $\fnn{h}{S_{\bbR^n}}{\bbR^{n - 1}}$
there is $x \in S_{\bbR^n}$ such that $h(x) = h(-x)$.
\rule{0pt}{0pt}\hfill\proofend

\begin{cor}\label{c3.6}
Let $X$ be an infinite dimensional normed space.
Then $S_X$ does not have a cover $\calS$ such that
$\calS$ has finite order, $\calS$ is locally finite,
and for every $S \in \calS$,
$S$ is a closed and convex and $0 \not\in S$.
\vspace{-2.3mm}
\end{cor}

{\bf Proof }
Let $k$ be the order of $\calS$.
Let $n$ be such that $\frac{n}{2} + 1 > k$, and $L$ be an
$n$-dimensional subspace of $X$.
Let
$\calF =
\setm{S \cap S_L}{S \in \calS}$.
$L$ and $\calF$ fulfill the conditions of Lemma \ref{l3.5}.
That is,
\newline
(Q1) Every member of $\calF$ is closed, and does not contain antipodal
points.
\newline
(Q2) $\bigcup_{F \in \calF} F = S_{L}$.
\newline
(Q3) $\calF$ is finite.
\newline
The argument that (Q3) holds is as follows.
$\calS$ is locally finite. This implies that $\calF$ is locally
finite.
But a locally finite family of subsets of a compact space must
be finite.

By Lemma \ref{l3.5}, the order of $\calF$ is $> \frac{n}{2} + 1 > k$.
This contradicts the fact that the order of $\calS$ is $k$.
\hfill\proofend

\begin{prop}\label{p3.6}
\num{a}
Let $X$ be a normed space, $Y$ be a separable subspace of $X$ and
$g \in \itGamma_2(X)$. Then there is a separable subspace $Z$ of $X$
such that $Y \subseteq Z$ and $g(Z) \subseteq Z$.

\num{b} \num{i} Let $X$ be a normed space, $Z$ be a subspace of $X$
and $g \in \itGamma_2(X)$. Suppose that $g(Z) \subseteq Z$.
Then $g \restriction Z \in \itGamma_2(Z)$.

\num{ii} The same holds for $\itGamma_1$.
\vspace{-2.3mm}
\end{prop}

{\bf Proof } (a) Suppose first that $g \in \whatGamma_1(X)$.
Then $\calS_1 \eqdf \setm{S^g_j \cap Y}{j \in J^g}$\break
is locally finite.
$Y$ is second countable, and in a second countable Hausdorff space
every locally finite family has cardinality $\leq \aleph_0$.
So $\abs{\calS_1} \leq \aleph_0$.

For every $S \in \calS_1$, let $j_S \in J^g$ be such that
$S^g_j \cap Y = S$ and let $J_1 = \setm{j_S}{S \in \calS_1}$
and $E_1 = [\kern1pt \bigcup_{j \in J'} E^g_j \kern1pt]$.
So  $\abs{J_1} \leq \aleph_0$ and hence
$\fs{dim}(E_1) \leq \aleph_0$.
Let $Y_1 = Y + E_1$.
Clearly, $Y \cup g(Y) \subseteq Y_1$.
Since $Y$ is separable and $\fs{dim}(E_1) \leq \aleph_0$,
$Y_1$ is separable.

Appling the same process to $Y_1$
instead of $Y$ one obtains a subspace $Y_2$ such that
($*$) $Y_1 \cup g(Y_1) \subseteq Y_2$ and $Y_2$ is separable.

Repeating this procedure countably many times we get a sequenc
of separable subspaces
\newline\centerline{
$Y = Y_0 \subseteq Y_1 \subseteq Y_2 \subseteq \ldots$
}
Such that for $i = 0,1,\ldots,$ \ $g(Y_i) \subseteq Y_{i + 1}$.
Let $Z = \bigcup_{i = 0}^{\infty} Y_i$.
Then $Y \subseteq Z$,\break
$Z$ is separable and $g(Z) \subseteq Z$.

Assume now that $g \in \itGamma_1(X)$.
Let $g = h_n \scirc \ldots \scirc h_1$,
where $h_1,\ldots,h_n \in \whatGamma_1(X)$.
Let $\setm{f_i}{i \in \bbN}$ be an enumeration of
$\fsetn{h_1}{h_n}$ such that for every $m \leq n$
$\setm{i}{f_i = h_m}$ is infinite.
Let $Z_0 = Y$ and define by induction a chain of separable subspaces
$Z_0 \subseteq Z_1 \subseteq \ldots$ such that for every $i$,
$f_i(Z_i) \subseteq Z_{i + 1}$.
Let $Z = \bigcup_{i = 0}^{\infty} Z_i$.
Then $Y \subseteq Z$, $g(Z) \subseteq Z$ and $Z$ is separable.

(b) Part (b) is Trivial.
\hfill\proofend

\kern2mm

\setcounter{equation}{0}

\noindent
{\bf Proof of Theorem \ref{t3.1}(a) }
We first prove that Theorem \ref{t3.1}(a) is true for separable
spaces. So let $\rpair{X}{\norm{.}}$
be an infinite dimensional separable normed space.

Assume to the contrary that $\itGamma_1(X)$
is not non-shrinking.
There is a norm $\3norm{.}$ on $X$ such that
$\3norm{.}$ is equivalent to $\norm{.}$,
and $\3norm{.}$ fulfills the conclusion of Lemma~\ref{l3.3}.
By Proposition \ref{p3.2}(d), $\itGamma_1(X)$ is not non-shrinking
with repect to $\rpair{X}{\3norm{.}}$.
We may thus assume
that $\norm{.}$ fulfills the conclusion of Lemma~\ref{l3.3}.

Let $g \in \itGamma_1(X)$, $r \in (0,1)$ and
$x_1,\ldots,x_m \in X$ be such that
\begin{equation}\label{e}
g(B_X) \subseteq \bigcup_{k = 1}^m (x_k + r B_X)
\end{equation}
Let
\newline\centerline{
$g = h_n \scirc h_{n - 1} \scirc \ldots \scirc h_1$,
where $h_i \in \whatGamma_1(X), \ i = 1,\ldots,n$.
}
Every separable metric space is second countable,
and in a second countable Hausdorff space every locally finite
family has cardinlity $\leq \aleph_0$.
So $\abs{J^{h_i}} \leq \aleph_0$.

We may assume that for every $i$, \,$\abs{J^{h_i}} = \aleph_0$.
This can be achieved by adding to
$\setm{\rpair{S^h_j}{E^h_j}}{j \in J^h}$
additional pairs of the form $\rpair{S}{\sngltn{0}}$.
So denote the set of pairs associated with $h_i$ by
$\setm{\rpair{S_j^i}{E_j^i}}{j \in \bbN}$.

Let $\sigma_1 = \setm{j \in \bbN}{0 \in S_j^1}$.
Since $\setm{S_j^1}{j \in \bbN}$ is separated,
$\sigma_1$ is finite. (In fact, $\abs{\sigma_1} \leq 1$).
Let $M_1$ be the following finite-dimensional subspace of~$X$.
\newline\centerline{
$M_1 =
[\kern1pt \fsetn{x_1}{x_m} \cup
\bigcup_{j \in \sigma_1} E_j^1 \kern1pt]$.
}
(The points $x_1,\dots,x_m$ were defined in (\ref{e})).
Let $L_1\subseteq X$ be a complement of $M_1$ in $X$ and
$\fnn{P_1}{X}{L_1}$ be the projection of $\rpair{L_1}{M_1}$.
(That is, $P_1 \restriction M_1 = 0$
and
$P_1 \restriction L_1 = \fs{Id}$).

By Lemma \ref{l3.4}(c), $\setm{P_1(S^2_j)}{j \in \bbN}$
has finite order.
So
\newline\centerline{
$\sigma_2 \eqdf \setm{j \in \bbN}{0 \in P_1 (S_j^2)}$
}
is finite.
Put $M_2 = [\kern1pt M_1 \cup \bigcup_{j \in \sigma_2} E_j^2 \kern1pt]$,
and let $L_2 \subseteq X$ be a complement of $M_2$ in $X$.
Let $\fnn{P_2}{X}{L_2}$
be the projection of $\rpair{L_2}{M_2}$.

Denote $M_0 = [\kern1pt \fsetn{x_1}{x_m} \kern1pt]$
and $P_0 = \fs{Id}_X$.
Proceeding in the above way\break
$n$ times we construct:
\newline
A chain of finite-dimensional subspaces of $X$
\newline\centerline{
$M_1 \subseteq M_2 \subseteq \ldots \subseteq M_n$.
}
A sequence of complements of the $M_i$'s  with respect to $X$
\newline\centerline{
$L_1 ,\; L_2,\ldots,L_n$.
}
The projections
\newline\centerline{
$\fnn{P_i}{X}{L_i}$ \ \ of\ \ \ $\rpair{L_i}{M_i}$,
$\ \, i = 1,\ldots,n$.
}
And finite sets
\newline\centerline{
$\sigma_i =
\setm{j \in \bbN}{0 \in P_{i - 1}(S_j^i)},\; i = 1,\ldots,n$
}
such that
\newline\centerline{
$M_i =
[\kern1pt M_{i - 1} \cup
\bigcup_{j \in \sigma_i} E_j^i \kern1pt],\; i = 1,\ldots,n$.
}
(In fact, $L_n$ and $P_n$ will not be used).

Choose $\varepsilon > 0$ such that $r(1 + \varepsilon) < 1$.
Recall that $\norm{.}$ has the property of Lemma \ref{l3.3}.
Let $L$ be a subspace of $X$
as guaranteed by Lemma \ref{l3.3} for $M_n$ and $\varepsilon$.
That is,
\newline
(P1) $L$ has finite codimension in $X$ and $L \cap M_n = \sngltn{0}$.
\newline
(P2) Let $\fnn{P}{M_n + L}{L}$ be the projection of $\rpair{L}{M_n}$.
Then $\norm{P} < 1 + \varepsilon$.
\newline
We may also assume that
\newline
(P3) $L \subseteq \bigcap_{i = 1}^{n - 1} L_i$.
\newline
This is so, since $L \cap \bigcap_{i = 1}^{n - 1} L_i$ too, fulfills
(P1) and (P2).
\newline
Let
\begin{equation}\label{ee}
\mbox{$
C = S_L \setminus
\bigcup_{i = 1}^n \bigcup_{j \not\in \sigma_i} P_{i - 1}(S_j^i).
$}
\end{equation}

We show that $C \neq \emptyset$. Suppose otherwise. So
$S_L \subseteq
\bigcup_{i = 1}^n \bigcup_{j \not\in \sigma_i} P_{i - 1}(S_j^i)$.
By Lemma \ref{l3.4}(a), for every $i,j$, $P_{i - 1}(S_j^i)$ is closed.
If $j \not\in \sigma_i$,
then $0 \not\in P_{i - 1}(S_j^i)$.
Hence since $P_{i - 1}(S_j^i)$ is convex,
it does not contain antipodal points.
We use the facts that for every
$i = 1,\ldots,n$, $\setm{S^i_j}{j \in \bbN}$ is separated and bounded.
By Lemma \ref{l3.4}(c),
$\setm{P_{i - 1}(S_j^i)}{j \in \bbN \setminus \sigma_i}$
has finite order.
So
$\calS \eqdf
\setm{P_{i - 1}(S_j^i)}{i = 1,\ldots,n,\ \,j \not\in \sigma_i}$
has finite order.
By Lemma \ref{l3.4}(b),\break
$\calS$ is locally finite.
The facts:
\newline
(i) $\calS$ covers $S_L$;
\newline
(ii) For every $S \in \calS$, $S$ is closed and convex
and $0 \not\in S$;
\newline
(iii) $\calS$ is locally finite;
\newline
(iv) $\calS$ has finite order;
\newline
contradict Corollary \ref{p3.6}.
So $C \not= \emptyset$.

Recall that $P$ is the projection of $\rpair{L}{M_n}$.
We claim that
\newline
($*$)\kern4pt
\indent
For every $x \in C$, $g(x) = x + v$ for some $v \in M_n$.

Let $x \in C$.

{\it Step 1:}
By (\ref{ee}), either $x \not\in \bigcup_{j \in \bbN} S^1_j$,
or $x \in \bigcup_{j \in \sigma_1} S^1_j$.
In either case,
\newline\centerline{
$h_1(x) = x + v_1 ,\;  v_1 \in M_1$.
}
>From the facts that $x \in L \subseteq L_1$ and $v_1 \in M_1$,
it follows that $x = P_1(x + v_1)$.

{\it Step 2:}
Assume by contradiction that
$x + v_1 \in \bigcup_{j \not\in \sigma_2} S_j^2$.
Then
\newline\centerline{
$x = P_1(x + v_1) \in \bigcup_{j \not\in \sigma_2} P_1(S_j^2)$,
}
contradicting the fact that $x \in C$.
So $x + v_1 \not\in \bigcup_{j \not\in \sigma_2} S_j^2$.
Hence
either $x + v_1 \not\in \bigcup_{j \in \bbN} S^1_j$,
or $x + v_1 \in \bigcup_{j \in \sigma_2} S^2_j$.
In either case,
$$h_2 \scirc h_1(x) = h_2(x + v_1) = x + v_1 + v_2 ,\; v_2 \in M_2.$$
Since $M_1 \subseteq M_2$, $v_1 + v_2 \in M_2$.
Also, $x \in L \subseteq L_2$.
So $P_2(x + v_1 + v_2) = x$.

In {\it Step $n$} of this argument one concludes that

\kern 1mm

\centerline{
$g(x) =
x + \sum_{i = 1}^n v_i , \; v_i \in M_i,\; i = 1,\dots,n$.
}

\kern 1mm

\noindent
For every $i \leq n$, $M_i \subseteq M_n$.
So $v = \sum_{i = 1}^n v_i \in M_n$.
That is, $g(x) = x + v$, where $v \in M_n$. So ($*$) holds.

It follows from ($*$) that
\newline
($*$$*$) \indent For every $x \in C$, $Pg(x) = x$.

We check that
\newline
($*$$*$$*$)\kern-4pt \indent For every $x \in C$,
$g(x) \in \bigcup_{k = 1}^m x_k + r B_{M_n + L}$.
\newline
It is given that
$g(x) \in \bigcup_{k = 1}^m x_k + r B_X$.
Write $g(x) = x_k + r u$, where $u \in B_X$.
By ($*$), $g(x) \in M_n + L$, and from the definition of $M_1$
follows that
\newline
$x_k \in M_1 \subseteq M_n + L$.
So $u = \frac{g(x) - x_k}{r} \in M_n + L$.
So $u \in B_{M_n + L}$.
\newline
By ($*$$*$$*$),
\newline\centerline{
$g(C) \subseteq \bigcup_{k = 1}^m (x_k + rB_{M_n + L})$.
}
By ($*$$*$),
\newline\centerline{
$C = P(g(C)) \subseteq \bigcup_{k = 1}^m P(x_k + rB_{M_n + L})$.
}
Since $x_k \in M_n$, it follows that
$P(x_k + rB_{M_n + L}) = r P(B_{M_n + L})$. Hence

\kern1mm

\centerline{
$C \subseteq \bigcup_{k = 1}^m r P(B_{M_n + L})$.
}

\kern1mm

\noindent
We now use the facts that $\norm{P} < 1 + \varepsilon$
and $(1 + \varepsilon) r < 1$. So

\kern1mm

\centerline{
$C \subseteq \bigcup_{k = 1}^m r P(B_{M_n + L}) \subseteq
\bigcup_{k = 1}^m r (1 + \varepsilon) \ncdot B_{M_n + L}$.
}

\kern1mm

\noindent
We know that $\emptyset \neq C \subseteq S_L$.
But all the points in
$\bigcup_{k = 1}^m r (1 + \varepsilon) \ncdot B_{M_n + L}$
have norm $< 1$. A contradiction.

We have shown that the claim of Theorem \ref{t3.1}(a) is true for
separable spaces.

Let $X$ be any normed space. Suppose by contradiction that
$\itGamma_1(X)$ is not non-shrinking.
Let $g \in \itGamma_1(X)$, $r \in (0,1)$ and
$x_1,\ldots,x_m \in X$ be such that
$g(B_X) \subseteq \bigcup_{k = 1}^m (x_k + r B_X)$.
Let $Z$ be a separable subspace of $X$ such that
$[\kern1pt \fsetn{x_1}{x_k} \kern1pt] \subseteq Z$
and $g(Z) \subseteq Z$. Such $Z$ exists by Proposition \ref{p3.6}.
Relying on the fact that $x_1,\ldots,x_m \in Z$, we conclude that
$g(B_Z) \subseteq \bigcup_{k = 1}^m (x_k + r B_Z)$.
Also, $g \restriction Z \in \itGamma_1(Z)$.
So $\itGamma_1(Z)$ is not non-shrinking.
This contradicts the first part of the proof.
So Part (a) of the theorem is proved.
\rule{0pt}{0pt}\hfill\proofend

\bigskip

The following additional facts are needed in the proof of
Theorem~\ref{t3.1}(b).

\begin{lemma}\label{pl1}
{\rm(\cite{F} Crollary 3, \cite{FL} Theorem 2.3)}
Assume that an infinite dimensional Banach space $E$
contains a nonempty bounded open subset
which is a weak $G_{\delta}$-set.
Then $E$ contains $c_0$ isomorphically.
\end{lemma}

\begin{lemma}\label{l3.8}
Let $L$ be an infinite dimensional Banach space that does not contain
$c_0$ isomorphically,
and $\calA$ be a locally finite family of
$w$-closed subsets of $L$ which do not contain $0$.
Then $S_L \setminus \bigcup_{A \in \calA} A \not= \emptyset$.
\end{lemma}

\noindent
{\bf Proof } Assume to the contrary that
\newline\centerline{
$S_L \subseteq \bigcup_{A \in \calA} A$.
}
Let $E \subseteq L$ be any separable closed infinite dimensional
subspace of $L$.
Then $\calB = \setm{A \cap E}{A \in \calA}$ is a locally finite family
of w-closed subsets of $E$
and $S_E \subseteq \bigcup_{A \in \calB} A$.

Since $E$ is second countable $\abs{\calB} \leq \aleph_0$.
Define
\newline\centerline{
$G = B_E \setminus  \bigcup_{A \in \calB} A$.
}
Clearly, $0 \in G$.
Since $E$ is separable it follows that
$B_E$ is a weak $G_{\delta}$ set,
and hence $G$ is a weak $G_{\delta}$ set too.
Clearly $G$ is bounded.
We check that $G$ is open. Let $x \in G$.
Since $S_E \subseteq \bigcup_{A \in \calB} A$,
$x \in \fs{int}(B_E)$.
Since $\calB$ is a locally finite family
and $x \not\in \bigcup_{A \in \calB} A$,
it follows that there is $r > 0$ such that
$(x + r B_E) \cap \bigcup_{A \in \calB} A = \emptyset$.
Put $\alpha = \min(d(x,S_E),r)$.
Then $x + \alpha B_E \subseteq G$.
By Lemma \ref{pl1}, $c_0$ is embeddable in $E$. A contradiction.
\hfill\proofend

\medskip

\begin{lemma}\label{l3.9}
Let $M \subseteq X$ be a finite-dimensional subspace of a Banach space
$X$, $L \subseteq X$ be a closed complement of $M$ in $X$
and $\fnn{P}{X}{L}$ be the projection of $\rpair{L}{M}$.
Then for every w-closed bounded subset $A \subseteq X$,
$P(A)$ is w-closed.
\end{lemma}

\noindent
{\bf Proof } Denote the weak limit of a net $\calN$
by $\rfs{w-lim}\, \calN$.
We prove that $P(A)$\break
is closed under convergent nets.
Let $\rpair{D}{\leq_D}$ be a directed poset and\break
$\setm{x_d}{d \in D} \subseteq P(A)$ be a net in $A$
such that $\setm{x_d}{d \in D}$ is w-convergent in $L$, and
let $x = \rfs{w-lim} \setm{x_d}{d \in D}$.
Let $z_d \in A$ be such that $P(z_d) = x_d$.
So $z_d = x_d + y_d$, where $y_d \in M$.
The set $\setm{y_d}{d \in D}$ is bounded, since it is the image
of the bounded set $\setm{z_d}{d \in D}$ under a bounded operator.
So its closure is compact in $M$,
and hence the net $\setm{y_d}{d \in D}$ has a convergent subnet.
Let $\setm{y_{d_c}}{c \in C}$ be a convergent subnet of
$\setm{y_d}{d \in D}$ and $y = \rfs{w-lim} \setm{y_{d_c}}{c \in C}$.
So $x = \rfs{w-lim} \setm{x_{d_c}}{c \in C}$.
It follows that
$x + y = \rfs{w-lim} \setm{x_{d_c} + y_{d_c}}{c \in C}$.
Since $x_{d_c} + y_{d_c} = z_{d_c} \in A$ and $A$ is w-closed,
$x + y \in A$.
Clearly, $x = P(x + y) \in P(A)$.
We have shown that $P(A)$ is w-closed.
\hfill\proofend

\setcounter{equation}{0}

{\bf Proof of Theorem \ref{t3.1}(b) }
We first prove that Theorem \ref{t3.1}(b) is true for separable
spaces. So let $\rpair{X}{\norm{.}}$
be an infinite dimensional separable normed space.

Assume to the contrary that $\itGamma_2(X)$
is not non-shrinking.
There is a norm $\3norm{.}$ on $X$ such that
$\3norm{.}$ is equivalent to $\norm{.}$,
and $\3norm{.}$ fulfills the conclusion of Lemma~\ref{l3.3}.
By Proposition \ref{p3.2}(c), $\itGamma_2(X)$ is not non-shrinking
with repect to $\rpair{X}{\3norm{.}}$.
We may thus assume
that $\norm{.}$ fulfills the conclusion of Lemma~\ref{l3.3}.

Let $g \in \itGamma_2(X)$, $r \in (0,1)$ and
$x_1,\ldots,x_m \in X$ be such that
\begin{equation}\label{lfe}
g(B_X) \subseteq \bigcup_{k = 1}^m (x_k + r B_X)
\end{equation}
Let
\newline\centerline{
$g = h_n \scirc h_{n - 1} \scirc \ldots \scirc h_1$,
where $h_i \in \whatGamma_2(X), \ i = 1,\ldots,n$.
}
Every separable metric space is second countable,
and in a second countable Hausdorff space every locally finite
family has cardinlity $\leq \aleph_0$.
So $\abs{J^{h_i}} \leq \aleph_0$.

We may assume that for every $i$, \,$\abs{J^{h_i}} = \aleph_0$.
This can be achieved by adding to
$\setm{\rpair{\whatS^h_j}{E^h_j}}{j \in J^h}$
additional pairs of the form $\rpair{\whatS}{\sngltn{0}}$.
So denote the set of pairs associated with $h_i$ by
$\setm{\rpair{\whatS_j^i}{E_j^i}}{j \in \bbN}$.

Let $\sigma_1 = \setm{j \in \bbN}{0 \in \whatS_j^1}$.
Since $\setm{\whatS_j^1}{j \in \bbN}$ is locally finite in $\overX$,
$\sigma_1$ is finite.
Let $M_1$ be the following finite-dimensional subspace of $X$.
\newline\centerline{
$M_1 =
[\kern1pt \fsetn{x_1}{x_m} \cup
\bigcup_{j \in \sigma_1} E_j^1 \kern1pt]$.
}
(The points $x_1,\dots,x_m$ were defined in (\ref{lfe})).
Let $L_1\subseteq \overX$ be a complement of $M_1$ in $\overX$ and
$\fnn{P_1}{\overX}{L_1}$ be the projection of $\rpair{L_1}{M_1}$.
(That is, $P_1 \restriction M_1 = 0$
and
$P_1 \restriction L_1 = \fs{Id}$).

By Lemma \ref{l3.4}(b), $\setm{P_1(\whatS_j^2)}{j \in \bbN}$
is locally finite.
So
\newline\centerline{
$\sigma_2 \eqdf \setm{j \in \bbN}{0 \in P_1 (\whatS_j^2)}$
}
is finite.
Put
$M_2 = [\kern1pt M_1 \cup \bigcup_{j \in \sigma_2} E_j^2 \kern1pt]$,
and let $L_2 \subseteq \overX$ be a complement of $M_2$ in $\overX$.
Let $\fnn{P_2}{X}{L_2}$
be the projection of $\rpair{L_2}{M_2}$.

Denote $M_0 = [\kern1pt \fsetn{x_1}{x_m} \kern1pt]$
and $P_0 = \fs{Id}_{\overX}$.
Proceeding in the above way\break
$n$ times we construct:
\newline
A chain of finite-dimensional subspaces of $X$
\newline\centerline{
$M_1 \subseteq M_2 \subseteq \ldots \subseteq M_n$.
}
A sequence of complements of the $M_i$'s  with respect to $\overX$
\newline\centerline{
$L_1 ,\; L_2,\ldots,L_n$.
}
The projections
\newline\centerline{
$\fnn{P_i}{\overX}{L_i}$ \ \ of\ \ \ $\rpair{L_i}{M_i}$,
$\ \, i = 1,\ldots,n$.
}
And finite sets
\newline\centerline{
$\sigma_i =
\setm{j \in \bbN}{0 \in P_{i - 1}(\whatS_j^i)},\; i = 1,\ldots,n$
}
such that
\newline\centerline{
$M_i =
[\kern1pt M_{i - 1} \cup
\bigcup_{j \in \sigma_i} E_j^i \kern1pt],\; i = 1,\ldots,n$.
}
(In fact, $L_n$ and $P_n$ will not be used).

Choose $\varepsilon > 0$ such that $r (1 + \varepsilon) < 1$.
Recall that $\norm{.}$ has the property of Lemma \ref{l3.3}.
Let $L$ be a subspace of $X$
as guaranteed by Lemma \ref{l3.3} for $M_n$ and $\varepsilon$.
That is,
\newline
(P1) $L$ has finite codimension in $X$ and $L \cap M_n = \sngltn{0}$.
\newline
(P2) Let $\fnn{P}{M_n + L}{L}$ be the projection of $\rpair{L}{M_n}$.
Then $\norm{P} < 1 + \varepsilon$.
\newline
We may also assume that
\newline
(P3) $L \subseteq \bigcap_{i = 1}^{n - 1} L_i$.
\newline
This is so, since $L \cap \bigcap_{i = 1}^{n - 1} L_i$ too, fulfills
(P1) and (P2).
Then
by the assumptions of Part (b),
$X$ has an infinite dimensional subspace $Y$ such that $c_0$ is not
isomorphically embeddable in $\overY$.
Let $Z = Y \cap L$. Since $L$ has finite codimension in $X$,
$Z$ is infinite dimensional. Clearly,
\newline
(i) $c_0$ is not isomorphic to a subspace of $\overZ$.

Let
\begin{equation}\label{fsee}
\mbox{$
C = S_{\oversZ} \setminus
\bigcup_{i = 1}^n \bigcup_{j \not\in \sigma_i} P_{i - 1}(\whatS_j^i).
$}
\end{equation}

We show that $C \neq \emptyset$.
Suppose otherwise.
So
\newline
(ii)
$S_{\oversZ} \subseteq
\bigcup_{i = 1}^n \bigcup_{j \not\in \sigma_i} P_{i - 1}(\whatS_j^i)$.

Recall that $\bigcup_{j \in \bbN} \whatS^i_j$ is bounded.
By Lemma \ref{l3.9}, for every $i,j$,
$P_{i - 1}(\whatS_j^i)$ is w-closed in $\overX$.
So
\newline
(iii)
$P_{i - 1}(\whatS_j^i) \cap \overZ$ is w-closed in $\overZ$.

By Lemma \ref{l3.4}(b), for every $i$,
$\setm{P_{i - 1}(\whatS_j^i)}{j \in \bbN \setminus \sigma_i}$
is locally finite.
So
\newline
(iv)
$\setm{P_{i - 1}(\whatS_j^i) \cap \overZ}
{i = 1,\ldots,n, \ \, j \in \bbN \setminus \sigma_i}$
is locally finite.
\newline
By the definition of the $\sigma_i$'s,
\newline
(v)
$0 \not\in P_{i - 1}(\whatS_j^i)$ whenever $j \not\in \sigma_i$.
\newline
(i) - (v) contradict Lemma \ref{l3.8}.
So $C \neq \emptyset$.

We check that $C$ is open in $S_{\overZ}$.
This is so, since $C$ is the complement in $S_{\overZ}$
of the union of the locally finite family of closed sets
\newline
$\setm{P_{i - 1}(\whatS_j^i) \cap \overZ}
{i = 1,\ldots,n, \ \, j \in \bbN \setminus \sigma_i}$.

Since $Z \cap S_{\oversZ}$ is dense in $S_{\oversZ}$,
$C \cap Z \neq \emptyset$. In particular, $C \cap L \neq \emptyset$.

Recall that $P$ is the projection of $\rpair{L}{M_n}$.
We claim that
\newline
($*$)\kern4pt
\indent
For every $x \in C \cap L$, $g(x) = x + v$ for some $v \in M_n$.

Let $x \in C \cap L$.

{\it Step 1:}
By (\ref{fsee}), either $x \not\in \bigcup_{j \in \bbN} \whatS^1_j$,
or $x \in \bigcup_{j \in \sigma_1} \whatS^1_j$.
In either case,
\newline\centerline{
$h_1(x) = x + v_1 ,\;  v_1 \in M_1$.
}
>From the facts that $x \in L \subseteq L_1$ and $v_1 \in M_1$,
it follows that $x = P_1(x + v_1)$.

{\it Step 2:}
Assume by contradiction that
$x + v_1 \in \bigcup_{j \not\in \sigma_2} \whatS_j^2$.
Then
\newline\centerline{
$x = P_1(x + v_1) \in \bigcup_{j \not\in \sigma_2} P_1(\whatS_j^2)$,
}
contradicting the fact that $x \in C$.
So $x + v_1 \not\in \bigcup_{j \not\in \sigma_2} \whatS_j^2$.
Hence
either $x + v_1 \not\in \bigcup_{j \in \bbN} \whatS^1_j$,
or $x + v_1 \in \bigcup_{j \in \sigma_2} \whatS^2_j$.
In either case,
$$h_2 \scirc h_1(x) = h_2(x + v_1) = x + v_1 + v_2 ,\; v_2 \in M_2.$$
Since $M_1 \subseteq M_2$, $v_1 + v_2 \in M_2$.
Recalling that $x \in L \subseteq L_2$, we conclude that
$P_2(x + v_1 + v_2) = x$.

In {\it Step $n$} of this argument one concludes that

\kern 1mm

\centerline{
$g(x) =
x + \sum_{i = 1}^n v_i , \; v_i \in M_i,\; i = 1,\dots,n$.
}

\kern 1mm

\noindent
For every $i \leq n$, $M_i \subseteq M_n$.
So $v = \sum_{i = 1}^n v_i \in M_n$.
That is, $g(x) = x + v$, where $v \in M_n$. So ($*$) holds.

It follows from ($*$) that
\newline
($*$$*$) \indent For every $x \in C \cap L$, $Pg(x) = x$.

We check that
\newline
($*$$*$$*$)\kern-4pt \indent For every $x \in C \cap L$,
$g(x) \in \bigcup_{k = 1}^m x_k + r B_{M_n + L}$.
\newline
It is given that
$g(x) \in \bigcup_{k = 1}^m x_k + r B_X$.
Write $g(x) = x_k + r u$, where $u \in B_X$.
By ($*$), $g(x) \in M_n + L$, and from the definition of $M_1$
follows that
\newline
$x_k \in M_1 \subseteq M_n + L$.
So $u = \frac{g(x) - x_k}{r} \in M_n + L$.
So $u \in B_{M_n + L}$.
\newline
By ($*$$*$$*$),
\newline\centerline{
$g(C \cap L) \subseteq \bigcup_{k = 1}^m (x_k + rB_{M_n + L})$.
}
By ($*$$*$),
\newline\centerline{
$C \cap L = P(g(C \cap L)) \subseteq
\bigcup_{k = 1}^m P(x_k + rB_{M_n + L})$.
}
Since $x_k \in M_n$, it follows that
$P(x_k + rB_{M_n + L}) = r P(B_{M_n + L})$. Hence

\kern1mm

\centerline{
$C \cap L \subseteq \bigcup_{k = 1}^m r P(B_{M_n + L})$.
}

\kern1mm

\noindent
We now use the facts that $\norm{P} < 1 + \varepsilon$
and $(1 + \varepsilon) r < 1$. So

\kern1mm

\centerline{
$C \cap L \subseteq \bigcup_{k = 1}^m r P(B_{M_n + L}) \subseteq
\bigcup_{k = 1}^m r (1 + \varepsilon) \ncdot B_{M_n + L}$.
}

\kern1mm

\noindent
We know that $\emptyset \neq C \cap L \subseteq S_L$.
But all the points in
$\bigcup_{k = 1}^m r (1 + \varepsilon) \ncdot B_{M_n + L}$
have norm $< 1$. A contradiction.

We have shown that the claim of Theorem \ref{t3.1}(b) is true for
separable spaces.

Let $X$ be any normed space. Suppose by contradiction that
$\itGamma_2(X)$ is not non-shrinking.
Let $g \in \itGamma_2(X)$, $r \in (0,1)$ and
$x_1,\ldots,x_m \in X$ be such that
$g(B_X) \subseteq \bigcup_{k = 1}^m (x_k + r B_X)$.
Let $Z$ be a separable subspace of $X$ such that
$[\kern1pt \fsetn{x_1}{x_k} \kern1pt] \subseteq Z$
and $g(Z) \subseteq Z$. Such $Z$ exists by Proposition \ref{p3.6}.
Relying on the fact that $x_1,\ldots,x_m \in Z$, we conclude that
$g(B_Z) \subseteq \bigcup_{k = 1}^m (x_k + r B_Z)$.
Also, $g \restriction Z \in \itGamma_2(Z)$.
So $\itGamma_2(Z)$ is not non-shrinking.
This contradicts the first part of the proof.
So Part (b) of the theorem is proved.
\rule{0pt}{0pt}\hfill\proofend

\bigskip

\begin{prop}\label{p3.10}
Theorem \ref{t3.1}(b) does not hold for $X = c_0$.
\end{prop}

\noindent
{\bf Proof }
Let $V$ be the set of all the sequence in $c_0$ which have
only finitely many nonzero coordinates and in which every coordinate
is $0$, $2$ or $-2$.
For $v \in V$ let $S_v = v + B_{c_0}$.

It is not difficult to check that
$\setm{S_v}{v \in V}$ is a locally finite family and that
$\bigcup_{v \in V} S_v = 3B_{c_0}$.
For every $x \in 3B_{c_0}$, let $v_x$ be such that
$x \in S_v$.
Define a map $\fnn{h}{c_0}{c_0}$ as follows.
If $x \in 3B_{c_0}$, then $h(x) = x - v_x$,
and if $x \in c_0 \setminus 3B_{c_0}$, then $h(x) = x$.
Clearly, $h \in \itGamma_2(c_0)$ and $h(3B_{c_0}) = B_{c_0}$.
\hfill\proofend

\begin{question}\label{q3.11}
\begin{rm}
Is $H(c_0) \cap \itGamma_2(c_0)$ non-shrinking?
\end{rm}
\end{question}

\newpage

\section{Faithfulness in normed spaces and in\\ metrizable locally
convex spaces}\label{s4}

This section deals with two faithfulness theorems: the first
concerns with open subsets of normed spaces, and the second
with open subsets of metrizable locally convex spaces.
Recall that $K_1$ is the faithful class of Theorem~\ref{t1.2}(a).

For every normed space $E$ we shall define a subgroup $G_E$ of $H(E)$
such that $G_E \subseteq \itGamma_1(E) \cap \fs{LLIP}(E)$,
and for every nonempty open subset $X \subseteq E$ we define
$G_X = \setm{g \nrestriction X}{g \in G_E \mbox{ and }
\fs{supp}(g) \subseteq E}$.
In the first theorem, Theorem~\ref{t4.3}, we prove that
$\rpair{X}{G_X} \in K_1$.
Also, since $G_E \subseteq \itGamma_1(E)$,
$\rpair{X}{G_X}$ does not have small sets.


Every metrizable topological vector space has a metric
which is invariant under $+$. We shall deal only with such metrics.
We next define the group $G_X$ mentioned above.

\begin{defn}\label{d4.1}
\begin{rm}
(a) Let $\rpair{E}{\tau}$ be a metrizable topological vector space
and $d$ be metric on $E$ whose topology is $\tau$ and which is
invariant under~$+$. Then $\rpair{E}{d}$ is called
a {\it metric vector space}.
In particular a {\it metric locally convex space} is a
metric vector space which is locally convex.


(b) Let $X,Y$ be metric spaces and $\fnn{f}{X}{Y}$ be $\onetoone$.
The function $f$ is {\it bilipschitz}
if $f$ and $f\inverse$ are
Lipschitz functions.
$f$ is {\it locally bilipschitz} if for every $x \in X$
there is $U \in \fs{Nbr}(x)$ such that $f \nrestriction U$ is
bilipschitz.


(c) For a metric space $X$, $A \subseteq X$ and $r > 0$, denote
\newline
$B_X(A,r) = \bigcup_{x \in A} B_X(x,r)$.
For a vector space $E$ and $x,y \in E$ denote
\newline
$[x,y] = \setm{tx + (1 - t)y}{t \in [0,1]}$.

%

(d) Let $E$ be a metric locally convex space and $h \in H(E)$.
We say that $h$ is a {\it basic homeomorphism}
if there are $x,v \in E$ and $r > 0$ such that
\newline
(1) $\fs{supp}(h) \subseteq B_E([x,x + v],r)$;
\newline
(2) 
For every $z \in E$, there is $\lambda(z) \in \bbR$ such that
$h(z) = z + \lambda(z) \ncdot v$.
\newline
(3) $h(x) = x + v$.
\newline
Denote $x = x^h$, $x + v = y^h$, $v = v^h$, $r = r^h$
and $B([x,x + v],r) = S^h$.\break
Note that $h = \fs{Id}$ is a basic homeomorphism with $v^h = 0$.
\vspace{-2.0mm}
\end{rm}
\end{defn}

\kern2pt

\noindent
{\bf Faithfulness in Normed spaces}

\begin{defn}\label{d4.2}
\begin{rm}
Let $E$ be a normed space.
We define $\whatG_E \subseteq H(E)$.
A homeomorphism $h$ belongs to $\whatG_E$ if there is a sequence
of basic homeomorphisms $\setm{h_i}{i \in \bbN}$ such that
\newline
(1) $\setm{S^{h_i}}{i \in \bbN}$ is separated.
\newline
(2) $h \nrestriction (E \setminus \bigcup_{i \in \bbN} S^{h_i}) =
\fs{Id}$,
and for every $i \in \bbN$,
\,$h \nrestriction S^{h_i} = h_i \nrestriction S^{h_i}$.
\newline
(3) $\fs{supp}(h)$ is a bounded set.
\newline
(4) For every $i,j \in \bbN$, \,$r^{h_i} = r^{h_j}$.
\newline
(5) $h$ is locally bilipschitz.
\newline
Let
$\whatG^{\fss{LIP}}_E =
\setm{h \in \whatG_E}{h \mbox{ is bilipschitz}}$.
Define $G_E$ and $G^{\fss{LIP}}_E$ to be the subgroups of $H(E)$
generated by $\whatG_E$ and $\whatG^{\fss{LIP}}_E$ respectively.
Suppose that $X \subseteq E$ is open and nonempty, and define
\newline\centerline{
$G_X =
\setm{g \nrestriction X}{g \in G_E \mbox{ and }
g \nrestriction (E \setminus X) = \fs{Id}}$
}
and
\newline\centerline{
$G^{\fss{LIP}}_X =
\setm{g \nrestriction X}{g \in G_E^{\fss{LIP}} \mbox{ and }
g \nrestriction (E \setminus X) = \fs{Id}}$.
}
\vspace{-1.0mm}
\end{rm}
\end{defn}

Note that $\whatG_E \subseteq \whatGamma_1(E)$ and
$(\whatG_E)\inverse = \whatG_E$.
Hence $G_E \subseteq \itGamma_1(E)$.

\begin{theorem}\label{t4.3}
\num{a} For every normed space $E$ and an open nonempty subset
$X \subseteq E$,
$\rpair{X}{G_X} \in K_1$.

\num{b} For every Banach space $E$
and an open nonempty subset $X \subseteq E$,
$\rpair{X}{G^{\fss{LIP}}_X} \in K_1$.
\end{theorem}

\noindent
{\bf Remark }
The definition of $K_1$ implies that $K_1$ is {\it closed upwards}.
That is, if $\rpair{X}{G} \in K_1$ and $G \subseteq H \subseteq H(X)$,
then $\rpair{X}{H} \in K_1$.
So the fact that $\rpair{X}{G^{\fss{LIP}}_X} \in K_1$ implies that
$\rpair{X}{G_X} \in K_1$.

We need the following proposition. It appears in \cite{RY} as
Proposition~2.14(c).

\begin{prop}\label{p4.4}
There is a function $M(\ell,t)$
increasing in $\ell$ and decreasing in $t$
such that for every normed
space~$E$, $x,y \in E$ and  $r > 0$,
there is
$h \in H(E)$ such that:
\newline
\num{1} $\fs{supp}(h) \subseteq B([x,y],r)$;
\newline
\num{2}
$h(z) - z \in \fs{span}(\sngltn{y - x})$,
for every $z \in B([x,y],r)$;
\newline
\num{3} $h(x) = y$;
\newline
\num{4} $h$ is
$M(\norm{x - y},r)$-bilpschitz;
\newline
\num{5} $h \rest B(x,\frac{2r}{3}) =
\bfs{tr}_{y - x} \rest B(x,\frac{2r}{3})$.
(Recall that $\bfs{tr}_v$ denotes the function\break
\phantom{\num{5} }\kern2pt$x \mapsto x + v$).
\vspace{-2.3mm}
\end{prop}

Note that in the above proposition Clauses (1) - (3) imply
that $h$ is a basic homeomorphism.

For a metric space $Y$, $x \in Y$ and $r > 0$ let $\dashB_Y(x,r)$
denote the open ball of $Y$ with center at $x$ and radius $r$.

\begin{prop}\label{p4.5}
Let $E$ be a normed space and $X \subseteq E$ be an open nonempty set.
Let $x \in X$ and $r > 0$ be such that
$\dashB_E(x,4r) \subseteq X$.
Then every infinite separated subset $A \subseteq \dashB_E(x,r)$
is dissectable with respect to $G^{\fss{LIP}}_X$.
\vspace{-2.3mm}
\end{prop}

\noindent
{\bf Proof }
If a set $C$ has a dissectable subset, then $C$ itself is dissectable.
So we may assume that $A$ is countable.
Let $A = \setm{x_n}{n \in \bbN}$.
There is $e > 0$ and a subsequence $\sngltn{x'_n}$ of $\sngltn{x_n}$
such that for every $\varepsilon > 0$ there is $n_0$ such that
for every $m,n > n_0$,
$e - \varepsilon < \norm{x'_m - x'_n} < e + \varepsilon$.
We may thus assume that for every distinct $m,n$,
$\frac{7e}{8} < \norm{x_m - x_n} < \frac{9e}{8}$.
Denote $x_1 = u$. We show that $\dashB_E(u,\frac{e}{2})$ dissects $A$.
By removing $x_1$ from $A$ we may assume that for every~$n$,
$\frac{7e}{8} < \norm{x_n - u} < \frac{9e}{8}$.

{\bf Claim 1 } For every $0 < s < \frac{e}{2}$
there are $h \in G^{\fss{LIP}}_E$ and $a > 0$ such that
\newline
(i) $\fs{supp}(h) \subseteq X$,
\newline
(ii) for every $n$,
$h(x_{2n}) \in \dashB_E(u,s)$
and $h \nrestriction \dashB_E(x_{2n + 1},a) = \fs{Id}$.
\newline
{\bf Proof } For simplicity assume that $u = 0$.
Let $y_n = \frac{s}{2} \ncdot \frac{x_n}{\norm{x_n}}$
and $L_n = [x_n,y_n]$.
We leave to the reader to check that there are $b,c > 0$
such that:
\newline
(1) \centerline{
$B_E(L_n,b) \subseteq X$,
}
\newline
(2) For every distinct $m,n$,
\newline\centerline{
$d(B_E(L_m,b),B_E(L_n,b)) > c$.
}
For every even $n \in \bbN$
let $h_n$ be as assured by Proposition~\ref{p4.4}.
That is,\break
\num{1} $\fs{supp}(h_n) \subseteq B_E(L_n,b)$;
\newline
\num{2} $h_n(x_n) = y_n$;
\newline
\num{3} $h_n$ is
$M(\norm{x_n - y_n},b)$-bilipschitz.

Note that because of the increasingness of $M(\ell,t)$ in $\ell$
there is $K$ such that for every $n \in 2 \bbN$,
$h_n$ is $K$-bilipschitz.
Let
\newline\centerline{
$h\  \ =\ \ %
\bigcup_{n \in 2 \bbN} h_n \nrestriction B_E(L_n,b)
\ \ \cup\ \ %
\fs{Id} \nrestriction (E \setminus \bigcup_{n \in 2 \bbN} B_E(L_n,b))$.
}
It is easy to see that $h$ is as required in Claim 1.

\kern1mm

We now show that $\dashB_E(x,\frac{e}{2})$ dissects $A$.
Let $U$ be a nonempty open subset of $\dashB_E(u,\frac{e}{2})$
and $y \in U$. There is $t > 0$ such that
$B_E([u,y],t) \subseteq B_E(u,\frac{e}{2})$.
By Proposition~\ref{p4.4}, there is $g \in G^{\fss{LIP}}_E$
such that $g(u) = y$ and $\fs{supp}(g) \subseteq B_E(u,\frac{e}{2})$.
Let $s > 0$ be such that $g(B_E(u,s)) \subseteq U$.
By Claim 1, there is $h \in G^{\fss{LIP}}_E$
such that $\fs{supp}(h) \subseteq X$,
for every $n \in 2 \bbN$, $h(x_n) \in \dashB(u,s)$
and for every $n \in 2 \bbN\kern-1pt + 1$,
there is $T \in \fs{Nbr}(x_n)$
such that $h \nrestriction T = \fs{Id}$.
Clearly, $\fs{supp}(g \scirc h) \subseteq X$.
Let $f = (g \scirc h) \nrestriction X$.
Then
\newline
(1) $f \in G^{\fss{LIP}}_X$;
\newline
(2) $f(x_n) \in U$, for every $n \in 2 \bbN$;
\newline
(3) If $n \in 2 \bbN\kern-1pt + 1$, then there is $T \in \fs{Nbr}(x_n)$
such that $f \nrestriction T = \fs{Id}$.
We have shown that $\dashB_E(x,\frac{e}{2})$ dissects $A$.
\hfill\proofend

Let $E$ be a normed space and $g \in \fs{LIP}(E)$.
Then $g$ can be extended uniquely
to a homeomorphism $\barg$ of $\overE$.
For simplicity, if $x \in \overE$ we denote $\barg(x)$ by $g(x)$.

\begin{prop}\label{p4.6}
Let $E$ be a normed space, $u \in E$
and $x \in \dashB_{\oversE}(u,r) \setminus E$.
Let $\emptyset \neq V \subseteq \dashB_{\oversE}(u,r)$ be open.
Then there is $h \in G^{\fss{LIP}}_E$ such that
\newline
$\fs{supp}(h) \subseteq \dashB_E(u,r)$ and $h(x) \in V$.
\vspace{-2.3mm}
\end{prop}

\noindent
{\bf Proof } Let $v \in V \cap E$.
Choose $w \in E$ sufficiently close to $x$ so that
{\thickmuskip=2.1mu \medmuskip=1.5mu \thinmuskip=1mu
$B_E([w,v],3 \norm{w - x}) \subseteq B_E(u,r)$ and
$\norm{w - x} < d(v,\overE \setminus V)$.
Denote $s = \norm{w - x}$.
}
\indent
We now use Proposition \ref{p4.4}. Let $h \in H(E)$ be such that
\newline
\num{1} $\fs{supp}(h) \subseteq B([w,v],3s)$;
\newline
\num{2}
$h(z) - z \in \fs{span}(\sngltn{v - w})$,
for every $z \in B([w,v],3s)$;
\newline
\num{3} $h(w) = v$;
\newline
\num{4} $h$ is
$M(\norm{w - v},3s)$-bilpschitz;
\newline
\num{5}
$h \rest B(w,2s) =
\bfs{tr}_{v - w} \rest B(w,2s)$.
\newline
Hence $h \in G^{\fss{LIP}}_E$ and $\fs{supp}(h) \subseteq B(u,r)$.
Also,
\newline
$h(x) = h(w) + (x - w) = v + (x - w)$.
So $\norm{h(x) - v} = \norm{x - w} < d(v,\overE \setminus V)$.
So $h(x) \in V$.
\hfill\proofend

Let $Y$ be a metric space and $\overY$ be its completion.
Let $U \subseteq Y$ be open. Denote
$\fs{cmpl}(U) = \setm{x \in \overE}
{\mbox{there is } V \in \fs{Nbr}_{\oversY}(x)
\mbox{ such that } V \cap Y \subseteq U}$.
Note that $\fs{cmpl}(U)$ is open in $\overE$ and that
$\fs{cmpl}(U) \cap E = U$.

\begin{prop}\label{p4.7}
Let $E$ be a normed space and $X \subseteq E$ be open.
Let\break
$\setm{x_n}{n \in \bbN} \subseteq X$
be a Cauchy sequence such that
$\lim_n x_n \in \fs{cmpl}(X) \setminus X$.
Then $\setm{x_n}{n \in \bbN}$
is dissectable with respect to $G_X$.
\vspace{-2.3mm}
\end{prop}

\noindent
{\bf Proof }
Denote $x = \lim_n x_n$.
Let $u \in X$ and $r > 0$ be such that
$\dashB_E(u,4r) \subseteq X$ and $\norm{x - u} < \frac{r}{4}$.
We may assume that the $x_n$'s are pairwise distinct and that
for every $n$, \,$\norm{x_n - x} < \frac{r}{4}$.
Let $v \in E$ be such that $\norm{v} = r$, $L = [x,x + v]$
and $L_n = [x_n,x_n + v]$.
Then $L \subseteq B_{\oversE}(u,2r) \setminus E$
and $L_n \subseteq B_E(u,2r)$.
There are $t > s > 0$ such that $B_E(L_1,t) \subseteq B_E(u,4r)$
and for every $n$, \,$B_E(L_n,s) \subseteq B_E(L_1,t)$,
(for example, take $t$ and $s$ to be $\frac{3r}{4}$ and $\frac{r}{4}$).
Let $\sngltn{s_n}$ be a sequence of positive numbers
such that for every $n \neq m$,
$B_E(L_n,s_n) \cap B_E(L_m,s_m) = \emptyset$
and $B_E(L_n,s_n) \subseteq B_E(L_1,t)$.
For every $n \in 2 \bbN$
let $h_n$ be as assured by Proposition~\ref{p4.4}
and such that $h_n(x_n) = x_n + v$
and $\fs{supp}(h_n) \subseteq B_E(L_n,s_n)$.
Let
\newline\centerline{
$h\  \ =\ \ %
\bigcup_{n \in 2 \bbN} h_n \nrestriction B_E(L_n,s_n)
\ \ \cup\ \ %
\fs{Id} \nrestriction
(E \setminus \bigcup_{n \in 2 \bbN} B_E(L_n,s_n))$.
}
Every accumulation point of
$\setm{B_E(L_n,s_n)}{n \in 2 \bbN}$ in $\overE$ belongs to $L$,
so since $L \cap E = \emptyset$,
$\setm{B_E(L_n,s_n)}{n \in 2 \bbN}$ is discrete in $E$.
It follows that $h \in H(E)$.

Since for every $n \in 2 \bbN$, $h_n$ is bilipschitz,
$h$ is locally bilipschitz.
Also, $\fs{supp}(h) \subseteq B_E([x_1,x_1 + v],t)$
and $h(x_1) = x_1 + v$.
Finally, for every $w \in E$, $h(w) - w \in \fs{span}(\sngltn{v})$.
So $h$ is a basic homeomorphism. Hence $h \in G_E$.

Note that $\fs{supp}(h) \subseteq B_E(L_1,t) \subseteq X$.
Also, for every $n \in 2 \bbN\kern-1pt + 1$,
$B_E(L_n,s_n) \in \fs{Nbr}(x_n)$
and $h \nrestriction B_E(L_n,s_n) = \fs{Id}$.

Let $w \in E$ be such that $\norm{(x + v) - w} < \frac{r}{4}$.
We show that $\dashB_E(w,\frac{r}{2})$
dissects $\setm{x_n}{n \in \bbN}$.

Clearly, $x + v \in \dashB_{\oversE}(w,\frac{r}{2})$,
$d(x,B_{\oversE}(w,\frac{r}{2}) > \frac{3r}{4}$
and $B_E(w,\frac{r}{2}) \subseteq X$.

Let $\emptyset \neq V \subseteq \dashB_E(w,\frac{r}{2})$ be open.
By Proposition~\ref{p4.6}, there is $g \in G_E$ such that
$g(x + v) \in V$ and $\fs{supp}(g) \subseteq \dashB_E(w,\frac{r}{2})$.
So $\fs{supp}(g) \subseteq X$.
It follows that $(g \scirc h) \nrestriction X \in G_X$.

Since $\lim_{n \in \bbN} h(x_{2n}) = x + v$
and $g(x + v) \in V$, for all but finitely many $n$'s,
\,$g \scirc h(x_{2n}) \in V$.

Let $n$ be odd. Then there is $T \in \fs{Nbr}(x_n)$ such that
$h \nrestriction T = \fs{Id}$.
Also, $\norm{x_n - x} < \frac{r}{4}$.
So $d(x_n,B_E(w,\frac{r}{2})) > 0$.
Hence $d(x_n,\fs{supp}(g)) > 0$.
It follows that there is $S \in \fs{Nbr}(x_n)$
such that $g \scirc h \mrestriction S = \fs{Id}$.

Let $f = (g \scirc h) \nrestriction X$. We have shown that $f \in G_X$,
and that the sets
$\setm{n}{f(x_n) \in V}$
and
$\setm{n}{\mbox{there is } T \in \fs{Nbr}(x_n) \mbox{ such that }
f \nrestriction T = \fs{Id}}$ are infinite.
So $\dashB_E(w,\frac{r}{2})$ dissects $\setm{x_n}{n \in \bbN}$.
\hfill\proofend

\kern2pt

\noindent
{\bf Proof of Theorem \ref{t4.3} }
We prove Parts (a) and (b) together.

Let $X$ be an open subset of a normed space $E$.
Then $X$ is regular and first countable.
That is, (P1) of Definition~\ref{d1.1}(f) holds.

By Proposition \ref{p4.4}, $G^{\fss{LIP}}_X$ is locally moving.
Hence every subgroup of $H(X)$ containing $G^{\fss{LIP}}_X$
is locally moving.
That is, (P2) holds for every $H$ containing $G^{\fss{LIP}}_X$.

Similarly, Proposition \ref{p4.4} implies that
$DF(X,G^{\fss{LIP}}_X) = X$. Hence the same is true for every
$H$ containing $G^{\fss{LIP}}_X$.
That is, (P3) holds for every $H$ containing $G^{\fss{LIP}}_X$.

We show that (P4) holds.
Assume that $E$ is a Banach space.
Let $x \in X$ and $r > 0$ be such that $\dashB(x,4r) \subseteq X$.
We show that $\dashB(x,r)$ is flexible with respect to
$G^{\fss{LIP}}_X$.
Let $A \subseteq \dashB(x,r)$ be an infinite set which is
discrete in $X$. So $A$ is discrete in $E$. Then $A$ contains
an infinite separated subset $A'$. By Proposition~\ref{p4.5},
$A'$ is disscetable with respect to $G^{\fss{LIP}}_X$.
So $A$ is disscetable with respect to $G^{\fss{LIP}}_X$,
and hence $\dashB(x,r)$ is flexible.
Since $\setm{\dashB(x,r)}{\dashB(x,4r) \subseteq X}$ is a cover of $X$,
it follows that $\rpair{X}{G^{\fss{LIP}}_X}$ satisfies (P4),
and the same is true for every $H$ containing $G^{\fss{LIP}}_X$.
We have proved Part (b) of Theorem~\ref{t4.3}.

Assume next that $E$ is a normed space.
Let $x \in X$ and $r > 0$ be such that $\dashB(x,4r) \subseteq X$.
We show that $\dashB(x,r)$ is flexible with respect to
$G_X$.
Let $A \subseteq \dashB(x,r)$ be an infinite set which is
discrete in $X$. Then $A$ is discrete in $E$.
Then either (i) $A$ contains an infinite separated subset $A'$,
or (ii) $A$ contains a Cauchy sequence $\setm{x_n}{n \in \bbN}$
converging to a point in $\overE \setminus E$.
If (i) happens, then by Proposition~\ref{p4.5}, $A'$ is dissectable
with respect to $G^{\fss{LIP}}_X$.
Hence $A$ is dissectable with respect to $G_X$.

If (ii) happens, then by Proposition~\ref{p4.7},
$\setm{x_n}{n \in \bbN}$ is dissectable
with respect to $G_X$.
Hence $A$ is dissectable with respect to $G_X$.
So $\dashB(x,r)$ is flexible with respect to $G_X$.
Hence $\rpair{X}{G_X}$ satisfies (P4).
We have proved Part (a) of Theorem~\ref{t4.3}.
\hfill\proofend

\kern3pt

\noindent
{\bf Faithfulness in metrizable locally convex spaces}

{\thickmuskip=2mu \medmuskip=1mu \thinmuskip=1mu
We now turn to general locally convex spaces.
The question~whether the class of locally convex topological
vector spaces is faithful is open. That is, it is unknown whether
$\setm{\rpair{X}{H(X)}}{X \mbox { is a locally convex space}}$
is faithful.
}

In \cite{LR} it was shown that if $K_2$
is the class of space-group pairs $\rpair{X}{G}$ such that
$X$ is an open subset of a normal locally convex space $E$,
and $E$ has a nonempty open set which intersects every line
in a bounded set. Then $K_2$ is faithful.

So $K_2$ includes spaces which are not first countable.
On the other hand, spaces which are a countable product of
normed spaces, and in particular, $\bbR^{\aleph_0}$ do not belong
to $K_2$.
These spaces do belong to the faithful class considered below.

Theorem D in the introduction states that
the class of all space-group
pairs $\rpair{X}{G}$
in which $X$ is an open subset of a locally convex metrizable
topological vector space, and $G$ is a group containing all locally
bi-uniformly continuous homeomorphisms of $X$ is faithful.
Below we define the group $P_X$ and prove the $\rpair {X}{P_X} \in K_1$.
This implies Theorem D.

\begin{defn}\label{d4.8}
\begin{rm}
Let $E$ be a metric locally convex space.

(a) Let $h \in H(E)$.
Suppose that there are $x_0,\ldots,x_n$, $k_1,\ldots,k_n$ and $r > 0$
such that $h = k_n \scirc \ldots \scirc k_1$,
and for every $i$, \,$k_i$ is a basic homeomorphism
with $x^{k_i} = x_{i - 1}$, $y^{k_i} = x_i$ and $r^{k_i} = r$.
Then $h$ is called a {\it polygonal homeomorphism}.
Denote $S^h = \bigcup_{i = 1}^n S^{k_i}$ and $r^h = r$.

(b)
We define $\whatP_E \subseteq H(X)$.
A member $h \in H(E)$ belongs to $\whatP_E$
if there is a set of polygonal homeomorphisms
$\setm{h_i}{i \in \bbN}$ such that
\newline
(1) $h \nrestriction (E \setminus \bigcup_{i \in \bbN} S^{h_i}) =
\fs{Id}$,
and for every $i \in \bbN$,
\,$h \nrestriction S^{h_i} = h_i \nrestriction S^{h_i}$.
\newline
(2) $h$ is locally bi-uniformly continuous.
\newline
(3) $\fs{supp}(h)$ is a bounded set.
\newline
(4) $\setm{S^{h_i}}{i \in \bbN}$ is discrete.
\newline
Let $P_E$ be the subgroup of $H(E)$ generated by $\whatP_E$
\newline
Let $X \subseteq E$ be open and nonempty. Define
\newline\centerline{
$P_X =
\setm{g \nrestriction X}{g \in P_E \mbox{ and }
g \nrestriction (E \setminus X) = \fs{Id}}$.
}
\vspace{-2.0mm}
\end{rm}
\end{defn}

\kern-08mm

\begin{theorem}\label{t4.9}
\num{a} Let $K_M$ be the class of all space-group
pairs $\rpair{X}{G}$ in which $X$ is a nonempty open subset of
a metrizable locally convex topological vector space $E$,
and $G \subseteq H(X)$ is a group containing all
locally bi-uniformly continuous homeomorphisms of~$X$.
\newline
Then $K_M$ is faithful.

\num{b} For every metric locally convex space $E$ and a nonempty open
subset $X \subseteq E$, \,$\rpair{X}{P_X} \in K_1$.

\num{c} Let
$K_P =
\setm{\rpair{X}{G}}{X \mbox{ is a nonempty open subset of a}
\mbox{ metric locally}\break \mbox{convex space and }
P_X \subseteq G \subseteq H(X)}$.
\newline
Then $K_P$ is faithful.

\end{theorem}

\noindent
{\bf Remarks } (a) In \ref{t4.9}(a) the ``uniform continuity'' is with
repect to the uniformity of the topological group $\rpair{E}{+}$.
However, for every metric $d$ on $E$, if $d$ is invariant under $+$
and $d$ induces on $E$ the original topology, then the uniformity of
$\rpair{E}{d}$ is equal to the uniformity of $\rpair{E}{+}$.

(b) Because $K_1$ is closed upwards, \ref{t4.9}(b) implies
\ref{t4.9}(c). Since $K_M \subseteq K_P$, \ref{t4.9}(c) implies
\ref{t4.9}(a).

(c) The proof of \ref{t4.9}(b) is analogous to the proof of
Theorem \ref{t4.3}(a).

\kern2mm

For a metric locally convex space $\rpair{E}{d}$ and
$x \in E$, denote $\norm{x} = d(0,x)$.

\begin{prop}\label{p4.10}
Let $\rpair{E}{d}$ be a metric locally convex space.
Let $x,y \in E$ and $r > 0$.
Then there is a bi-uniformly continuous basic homeomorphism $h$
such that $h(x) = y$ and $\fs{supp}(h) \subseteq B([x,y],r)$.

\vspace{-2.3mm}
\end{prop}

\noindent
{\bf Proof }
We may assume that $x = 0$.
Let $L = \fs{span}(\sngltn{y})$.
Since $E$ is locally convex, there is a complement $F$ of $L$
such that the projection $P$ of $\rpair{F}{L}$ is continuous.
In a metric vector space with an invariant metric every continuous
linear operator is uniformly continuous. In particular, $P$ is
uniformly continuous.

By the triangle inequality and the fact that the metric $d$ is
$+$-invariant,
$B_L(0,\frac{r}{2}) + B_F(0,\frac{r}{2}) \subseteq B_E(0,r)$.
So
$B_L([0,y],\frac{r}{2}) + B_F(0,\frac{r}{2}) \subseteq
B_E([0,y],r)$.
Let $p = \fs{min}(\frac{r}{2},\frac{\norm{y}}{2})$.
Choose $y_1 \in [0,y]$ such that $\norm{y_1} = p$, and let
$q$ be such that $q y = y_1$.
We define a function $\fnn{g}{[0,p] \times \bbR}{\bbR}$ as follows.
Let $t \in [0,p]$. We first define a function
$\fnn{g_t}{\bbR}{\bbR}$.
The function
$g_t$ is the unique piecewise linear function satisfying:
\newline\indent
(1) The breakpoints of $g_t$ are $-q,0$ and $1 + q$,
\newline\indent
(2) $g_t(s) = 0$ for every $s \in [-\infty,-q] \cup [1 + q,\infty]$,
\newline\indent
(3) $g_t(0) = \frac{p - t}{p}$.
\newline
Define $g$ by $g(t,s) = g_t(s)$.
Clearly, $g$ is uniformly continuous.

Let $z \in E$. Suppose that $z = \lambda y + v$,
where $v \in F$.
So $\lambda y = P(z)$.
Denote $\varphi(z) = \lambda$ and $v = Q(z)$.
Since $Q = \fs{Id} - P$, $Q$ is uniformly continuous.
$\varphi$ is a continuous homomorphism from $\rpair{E}{+}$ to
$\rpair{\bbR}{+}$. So $\varphi$ is uniformly continuous.
Define
\newline\centerline{
$h(z) =
\left\{
\begin{array}{ll}
g_{\snorm{v}}(\lambda) \ncdot y + v  &\indent\norm{v} \leq p \\
z                                    &\indent\norm{v} > p
\end{array}
\right.
$
}

\kern2mm

\noindent
Note that if $\norm{v} = p$, then $h(z) = z$.
Now, if $\norm{v} \leq p$, then\break
$h(z) = g(\norm{Q(v)},\varphi(z)) \ncdot y + Q(z)$.
Hence $h$ is uniformly continuous.
Clearly,

\kern2mm

\centerline{
$h\inverse(z) =
\left\{
\begin{array}{ll}
(g_{\lnorm{Q(z)}})\inverse(\varphi(z)) \ncdot y + Q(z)
&\indent\norm{Q(z)} \leq p \\
z                                    &\indent\norm{Q(z)} \geq p
\end{array}
\right.
$
}

\kern2mm

Since the function $\rpair{t}{s} \mapsto (g_t)\inverse(s)$
is uniformly continuous, $h\inverse$ is uniformly continuous.

Let $z \in E$. Suppose that $z = \lambda y + v$, where $v \in F$.
If $\norm{v} \geq p$, then $h(z) = z$,
and if $\lambda \not\in [-q,1 + q]$, then $h(z) = z$.
It follows that
\newline
$\fs{supp}(h) \subseteq B_L([0,y],p) + B_F(0,p) \subseteq
B_E([0,y],r)$.
So $h$ is as required.
\hfill\proofend

\begin{defn}\label{d4.11}
\begin{rm}
Let $E$ be a metrizable infinite dimensional locally convex space.
Denote the completion of $E$ by $\overE$.
A nonempty open subset
$U \subseteq E$ is called
{\it polygonally flexible},
if the following holds.
For every infinite $A \subseteq U$, if $A$ is discrete in $\overE$,
then there are an infinite set $B \subseteq A$
and $x \in U$
such that ${\rm P}(x,B)$ holds, where
${\rm P}(x,B)$ is the following statement.

For every $W \in \fs{Nbr}(x)$,
there are a set $\setm{y_b}{b \in B} \subseteq W$
and a family of polygonal lines $\setm{L_b}{b \in B}$
such that:
\begin{itemize}{}
\setlength{\parskip}{1.2mm}
\setlength{\itemsep}{0.0mm}
\item[(1)] For every $b \in B$, $L_b \subseteq \fs{cl}(U)$
and the endpoints of $L_b$ are $b$ and $y_b$.
\item[(2)] $\setm{L_b}{b \in B}$ is discrete in $\overE$.
\end{itemize}
\end{rm}
\end{defn}

\begin{prop}\label{p4.12}
\num{a} The separable Hilbert space $\ell_2$ is polygonally
flexible.
Moreover, for every infinite discrete $A \subseteq \ell_2$
there is an infinite $B \subseteq A$ and a nonempty open set $D$
such that for every $x \in D$, ${\rm P}(x,B)$ holds.

\num{b}
Let $E$ be an infinite dimensional metric locally convex space
and\break
$U \subseteq E$ be a convex open set.
Then $U$ is polygonaly flexible.
\vspace{-2.3mm}
\end{prop}

\noindent
{\bf Proof }
(a) The proof of Part (a) is easy and is left to the reader.

(b) Let $A \subseteq U$ be a countably infinite set discrete in
$\overE$.
Choose a separable infinite dimensional closed subspace $E_1$ of $E$
such that $E_1$ contains the linear span of $A$.
Let $U_1 = U \cap E_1$ and $F_1 = \fs{cl}_{\oversE_1}(U_1)$.

Clearly, $F_1$ is the closure of an open convex subset of a
separable Fr\'echet space.
By \cite{BP} Corollary 6.1 p.191,
such a set is homeomorphic to $\ell_2$.
Choose $\iso{h}{F_1}{\ell_2}$.
Let $A' = h(A)$. Then $A'$ is discrete in $\ell_2$.
So there is an infinite $B' \subseteq A'$ and a ball $D'$ of
$\ell_2$ such that for every $x' \in D'$,
${\rm P}(x',B')$ holds.
Let $x' \in D' \cap h(U_1)$. Let $x = h\inverse(x')$
and $B = h\inverse(B')$.
We show that ${\rm P}(x,B)$ holds.
Let $W \in \fs{Nbr}_E(x)$. Denote $W_1 = \fs{cmpl}(W \cap E_1)$.
Then $h(W_1) \in \fs{Nbr}_{\ell_2}(x')$.
Let $\setm{L'_{b'}}{b' \in B'}$ be as assured by ${\rm P}(x',B')$
for $h(W)$,
and $M_b = h\inverse(L'_{h(b)})$.
Replace $M_b$ by polygonal line $L_b$ which is sufficiently close
to $M_b$ and whose vertices are in $E$.
Then $L_b \subseteq E$.
Also, $\setm{L_b}{b \in B}$ is discrete in $F_1$.
Since $F_1$ is closed in $\overE$, $\setm{L_b}{b \in B}$ is discrete
in $\overE$.
\hfill\proofend

For $x \in E$ and $r > 0$ let $\dashB_E(x,r)$ denote the open ball
with center at $x$ and radius $r$.

\begin{prop}\label{p4.13}
Let $X$ be an open set
in an infinite dimensional metric locally convex space $E$.
Suppose that $U$ is an open nonempty bounded convex subset of $E$
such that $\fs{cl}_E(U) \subseteq X$.
Let $A \subseteq U$ be infinite,
and suppose that $A$ is discrete in $\overE$.
Then $A$ is dissectable in $\rpair{X}{P_X}$.
\vspace{-2.3mm}
\end{prop}

\noindent
{\bf Proof }
%
Let $B \subseteq A$ and $x \in U$ be such that ${\rm P}(x,B)$ holds.
We may assume that $\abs{B} = \aleph_0$ and that $x \not\in B$.

{\bf Claim 1 } For every $V \in \fs{Nbr}(x)$ and $C \subseteq B$
there is $h \in P_X$ such that for every $c \in C$, $h(c) \in V$,
and for every $b \in B \setminus C$ there is $W \in \fs{Nbr}(b)$
such that $h \nrestriction W = \fs{Id}$.
\newline
{\bf Proof }
For every $b \in B$ let $L_b \subseteq \fs{cl}_E(U)$ be a polygonal
line
and $y_b \in V$ be such that (i) the endpoints of $L_b$ are $b$ and
$y_b$, (ii) $\setm{L_b}{b \in B}$ is discrete in $\overE$.
Let
$d_b =
d(L_b,(E \setminus X) \cup
\bigcup \setm{L_{b'}}{b' \in B \setminus \sngltn{b}}$.
So $d_b > 0$.
Let $\setm{r_b}{b \in B}$ be a sequence of positive numbers converging
to $0$ such that $r_b < \frac{d_b}{2}$.
There is a bi-uniformly continuous polygonal homeomorphism $h_b$
such that $h_b(b) = y_b$ and $\fs{supp}(h_b) \subseteq B(L_b,r_b)$.
The homeomorphism $h_b$ is obtained by applying Proposition~\ref{p4.10}
to every edge of $L_b$.
So for every distinct $b,b' \in B$,
$\fs{supp}(h_b) \cap \fs{supp}(h_{b'}) = \emptyset$
and $\fs{supp}(h_b) \subseteq X$.
Also, since $\setm{r_b}{b \in B}$ converges to $0$,
$\setm{B(L_b,r_b)}{b \in B}$ is discrete in $\overE$.
Let
$h =
\fs{Id} \nrestriction (X \setminus \bigcup_{c \in C} B(L_c,r_c)) \cup
\bigcup_{c \in C} h_c$.
Then $h$ is as required.

{\bf Claim 2 } Let $d = \min(d(x,E \setminus U),\frac{d(x,B)}{2})$
and $T$ be an open convex set contained in $\dashB(x,d)$
and containing $x$.
Then $T$ dissects $B$.
\newline
{\bf Proof }
Let $W \subseteq T$ be open and nonempty.
Let $y \in W$. So $[x,y] \subseteq U$. Hence there is $r > 0$ such that
$B([x,y],r) \subseteq \dashB(x,d)$.
There is a bi-uniformly continuous basic homeomorphism $g'$
such that $\fs{supp}(g') \subseteq B([x,y],r)$ and $g'(x) = y$.
This follows from Proposition \ref{p4.10}.
Let $g = g' \nrestriction X$. Then $g \in P_X$.

Denote $V = g\inverse(W)$. Then $V \in \fs{Nbr}(x)$.
Let $C \subseteq B$ be such that $C$ and $B \setminus C$ are infinite.
By Claim~1, there is $h \in P_X$ such that
$h(C) \subseteq V$ and for every $b \in B \setminus C$,
there is $S \in \fs{Nbr}(b)$ such that $h \mrestriction S = \fs{Id}$.

It follows that for every $c \in C$, $g \scirc h(c) \in W$.
Also, $g \nrestriction (X \setminus B([x,y],r)) = \fs{Id}$
and $X \setminus B([x,y],r) \in \fs{Nbr}(b)$ for every $b \in B$.
So for every $b \in B \setminus C$ there is $S \in \fs{Nbr}(b)$
such that $g \scirc h \mrestriction S = \fs{Id}$.
This shows that $T$ dissects $B$.
This proves Claim 2.
Since $B \subseteq A$, $T$ dissects $A$ in $\rpair{X}{P_X}$.
\hfill\proofend

\kern1mm

If $g \in H(E)$ and $g$ is bi-uniformly continuous,
then $g$ has a unique extension to $\overE$.
So $g(x)$ is defined for any $x \in \overE$.

\begin{prop}\label{p4.14}
Let $E$ be a metric locally convex space.
Let $x \in \overE \setminus E$. Let $U$ be an open convex subset of
$E$ such that $x \in \fs{cmpl}(U)$.
Then for every nonempty open subset $V \subseteq U$
there is a bi-uniformly continuous basic homeomorphism $g$ of $E$
such that $\fs{supp}(g) \subseteq U$ and $g(x) \in \fs{cmpl}(V)$.
\vspace{-2.3mm}
\end{prop}

\noindent
{\bf Proof }
Let $z \in \fs{cmpl}(V)$ be such that $y \eqdf z - x \in E$.
Let $r > 0$ be such that
$B_{\oversE}([x,z],2r) \subseteq \fs{cmpl}(U)$.
So if $x' \in \overE$ is such that $\norm{x' - x} < r$,
then $B_{\oversE}([x',x' + y],r) \subseteq \fs{cmpl}(U)$.
By Proposition \ref{p4.10},
there is a bi-uniformly continuous basic homeomorphism $h$ of $E$
such that $\fs{supp}(h) \subseteq B([0,y],r)$ and $h(0) = y$.
Let $\barh$ be the extension of $h$ to $\barE$.
Let $e > 0$ be such that $B_{\oversE}(z,2e) \subseteq \fs{cmpl}(V)$.
There is $\varepsilon > 0$
such that for every $u \in \barE$:
if $\norm{u} < \varepsilon$, then
$\norm{\barh(u) - y} < e$.
Let $x' \in E$
be such that $\norm{x' - x} < \min(r,\varepsilon,e)$.
Let $g = \bfs{tr}_{x'} \scirc h \scirc \bfs{tr}_{x'}\inverse$.
So $g$ is a bi-uniformly continuous basic homeomorphism.
{\thickmuskip=3.5mu \medmuskip=2mu \thinmuskip=1mu
Clearly,
$\fs{supp}(g) \subseteq B([x',x' + y],r) \subseteq U$.
Since $\norm{x' - x} < \varepsilon$,
$\norm{g(x) - g(x')} < e$.
}
Also,
$g(x') = x' + z - x$.
So
$\norm{g(x) - x' - z + x} < e$.
Hence
\newline\centerline{
$\norm{g(x) - z} \leq \norm{g(x) - x' - z + x} + \norm{x' - x} < 2e$.
}
So $g(x) \in \fs{cmpl}(V)$.
\hfill\proofend

\begin{prop}\label{p4.15}
Let $E$ be a metric locally convex space.
Let $X \subseteq E$ be open,
and $\setm{x_n}{n \in \bbN} \subseteq X$
be a Cauchy sequence such that
\newline
$\lim_n x_n \in \fs{cmpl}(X) \setminus X$.
Then $\setm{x_n}{n \in \bbN}$
is dissectable with respect to $P_X$.
\vspace{-2.3mm}
\end{prop}

\kern-2mm

\noindent
{\bf Proof }
The proof is analogous to the proof of Proposition \ref{p4.7}.
It relies on Propositions \ref{p4.10} and \ref{p4.14} in the same way
that the proof of \ref{p4.7} relies on Propositions \ref{p4.4} and
\ref{p4.6}. \hfill\proofend

\kern 2mm

\noindent
{\bf Proof of Theorem \ref{t4.9} } Part (c) of \ref{t4.9} follows
from Part (b) and Part (a) follows from Part (c).
We prove Part (b).

Let $X$ be a nonempty open subset of a metric locally convex space
$E$.
Then $X$ is regular and first countable.
That is, (P1) of Definition~\ref{d1.1}(f) holds.
By Proposition \ref{p4.10}, $P_X$ is locally moving.
That is, (P2) holds for $P_X$.
Proposition \ref{p4.10} also implies that $DF_{X,P_X}(x,y)$ holds
for any distinct $x,y \in X$.
Indeed, choose $r > 0$ such that $y \not\in B_E(x,r) \subseteq X$.
Then by Proposition \ref{p4.10}, for every $z \in B^E(x,r)$
there is $g \in P_X$ such that $g(x) = z$,
and $g \nrestriction U = \fs{Id}$ for some $U \in \fs{Nbr}(y)$.
Hence (P3) holds.

We show that (P4) holds.
Let $U$ be a bounded convex open set such that
$\fs{cl}_{\oversE}(U) \subseteq \fs{cmpl}(X)$.
Let $A \subseteq U$ be an infinite set,
and assume that $A$ is discrete in $X$.
Either $A$ is discrete in $\overE$ or $A$ contains a $\onetoone$
Cauchy sequence converging to a member of $\overE - E$.
If $A$ is discrete in $\overE$, then by Proposition~\ref{p4.13},
$A$ is dissectable in $\rpair{X}{P_X}$.

Suppose that $B = \setm{x_n}{n \in \bbN} \subseteq A$ is
a $\onetoone$ Cauchy sequence converging to $x$ and
$x \in \overE - E$.  Then $x \in \fs{cmpl}(X) - X$.
By Proposition~\ref{p4.15}, $B$ is dissectable in $\rpair{X}{P_X}$.
We have shown that $U$ is flexible in $\rpair{X}{P_X}$.
So $X$ has an open cover consisting of flexilble sets.
That is, (P4) holds.
\hfill\proofend

\end{document}